\documentclass[12pt]{amsart}

\usepackage{amsmath,amssymb}
\usepackage{amsfonts}
\usepackage{graphicx}
\usepackage{epsfig}
\input amssym.def

\theoremstyle{plain} 
\newtheorem{proposition}{Proposition}[section]
\newtheorem{lemme}[proposition]{Lemme}
\newtheorem{remarque}[proposition]{Remarque}
\newtheorem{definition}[proposition]{D\'efinition}  

\newtheorem{corollaire}[proposition]{Corollaire}
\newtheorem*{corollaire }{Corollaire}

\newtheorem*{thmA}{Th\'eor\`eme A}
\newtheorem*{propA}{Proposition A}
\newtheorem*{thmA'}{Th\'eor\`eme A'}

\newtheorem*{thm }{Th\'eor\`eme}
\newtheorem*{non errance}{Conjecture de Non Errance}
\newtheorem*{propriete de point fixe}{Propri\'et\'e de Point fixe}

\newcommand{\diam}{\mbox{\rm diam}}

\newcommand{\can}{{\mbox{\rm can}}}

\renewcommand{\mod}{\, \, \mbox{\rm mod}}

\def\ov{\overline}
\def\cal{\mathcal}

\def\preuve{\par\noindent {\it Preuve.} }
\def\pp{{{\Bbb P}(\C_p)}}
\def\f{{\ov{\F}_p}}
\def\pf{{{\Bbb P}(\ov{\F}_p)}}

\def\e{\varepsilon}

\def\C{\Bbb C}

\def\F{\Bbb F}
\def\H{\Bbb H}

\def\Q{\Bbb Q}
\def\R{\Bbb R}

\def\Z{\Bbb Z}

\def\ca{{\cal A}}

\def\ce{{\cal E}}

\def\cm{{\frak m}}
\def\co{{\cal O}}
\def\cp{{\cal P}}
\def\cq{{\cal Q}}
\def\cs{{\cal S}}

\def\ct{{\cal T}}

\def\hB{\widehat{B}}

\def\hX{\widehat{X}}

\def\tB{\widetilde{B}}

\def\tr{\widetilde{R}}

\def\ts{{\widetilde{\cs}}}

\begin{document}
\title[Sur la structure des ensembles de Fatou $p$-adiques.]{Sur la structure des ensembles de Fatou $p$-adiques.}

\author[J. Rivera-Letelier]{Juan Rivera-Letelier}
\address{J. Rivera-Leteleir \\
         Mathematics Department\\
         SUNY at Stony Brook \\
         Stony Brook, NY 11794-3660}
\email {rivera@math.sunysb.edu}
\date{17 Janvier, 2002.}

\maketitle

\setcounter{tocdepth}{1}
\tableofcontents

\section{Introduction.}
	Fixons un nombre premier $p > 1$ et soient $\Q_p$ le corps des nombres $p$-adiques et $\C_p$ la compl\'etion d'une cloture agl\'ebrique de $\Q_p$.
	
	L'{\it ensemble de Fatou} $F(R) \subset \pp$ d'une fonction rationnelle $R \in \C_p(z)$ de degr\'e au moins deux, est par d\'efinition l'ensemble des points ayant un voisinage $X \subset \pp$ satisfaisant
$$
	\cup_{n \ge 0} R^n(X) \mbox{ omet au moins $3$ points de $\pp$}~;
$$
voir~\cite{Hs}, \cite{Be non rec}, \cite{Be hyp}, \cite{Be comp} et~\cite{these}.
	L'ensemble de Fatou est ouvert par d\'efinition et il est aussi dense dans $\pp$, voir~\cite{these}.

	D'apr\`es un th\'eor\`eme de L.C.~Hsia, tout point de l'ensemble de Fatou a un voisinage o\`u la famille $\{ R^n \}_{n \ge 1}$ est uniform\'ement lipschitzienne par rapport \`a la distance chordale sur $\pp$, voir~\cite{Hs}.
	
	Dans le cas complexe on peut d\'efinir l'ensemble de Fatou de la m\^eme fa\c{c}on.
	Pour \'etudier l'ensemble de Fatou dans le cas complexe on \'etudie ces {\it composantes connexes}.
	Chaque composante connexe est associ\'ee \`a un comportement dynamique bien determin\'e~: soit attractif (bassins d'attraction), soit quasi-p\'eriodique (disques de Siegel et anneaux de Herman)~; voir~\cite{Su}, \cite{Mi} ou~\cite{Bea}.

	Dans le cas $p$-adique on ne peut pas faire la m\^eme analyse (au moins directement), car la droite $\pp$ est totalement disconnexe (toute composante connexe de $\pp$ consiste d'un point).
D'autre part, il y a des fonctions rationnelles dont leur ensemble de Fatou est \'egal \`a~$\pp$ et qui cependant ont au m\^eme temps des comportements attractifs {\it et} quasi-p\'eriodiques.

	Dans ce papier on introduit une notion convenable de {\it composantes de l'ensemble de Fatou} (D\'efinition~\ref{composantes Fatou}) et on donne une caract\'erisation des composantes {\it p\'eriodiques} (Th\'eor\`eme~A), qui est analogue au cas complexe~; voir~\cite{Su}, \cite{Mi} ou~\cite{Bea}.

\subsection{Composantes de l'ensemble de Fatou.}
	 Comme $\pp$ est totalement disconnexe on remplace la notion de ``composante connexe'' par celle de ``{\it composante analytique}'', d\'efinit \`a continuation.
	Rappelons qu'un affino\"{\i}de ouvert est une intersection finie non vide de boules ouvertes.
\begin{definition}[Composante Analytique]\label{composante analytique}
	Soit $U$ une partie ouverte de $\pp$ et soit $x \in U$.
	Alors la {\bf composante analytique} de $U$ contenant $x$ est l'union de tous les {\rm affino\"{\i}des ouverts} $X \subset U$ contenant $x$.
\end{definition}
	Cette d\'efinition est plus restreint que celle introduite dans~\cite{Be comp}, voir aussi~\cite{these}.\footnote{Ici on consid\`ere {\it affino\"{\i}des ouverts} au lieu d'affino\"{\i}des ferm\'es.
	S'il on prend la d\'efinition comme dans~\cite{Be comp} ou~\cite{these} (avec affino\"{\i}des ferm\'es) l'\'enonc\'e correspondant du Th\'eor\`eme~A n'est pas vraie.
}

	Pour traiter le cas des fonctions rationnelles dont l'ensemble de Fatou est \'egal \`a $\pp$, on consid\`ere la d\'efinition suivante.
\begin{definition}[Composantes de l'ensemble de Fatou]\label{composantes Fatou}
	Soit $R \in \C_p(z)$ une fonction rationnelle de degr\'e au moins deux et soit $x \in F(R)$.
	Alors la {\bf composante de $F(R)$ contenant} $x$, est l'union de tous les affino\"{\i}des ouverts $X \subset F(R)$ contenant $x$ et tel que
$$
	\cup_{n \ge 0} R^n(X) \mbox{ omet au moins $3$ points de $\pp$.}
$$
\end{definition}
\noindent
	(On montre que deux composantes de l'ensemble de Fatou qui s'intersectent, co\"{\i}ncident.
Ceci n'est pas une cons\'equence imm\'ediate de la d\'efinition).

	{\it Dans le cas o\`u $F(R) \neq \pp$, les composantes de $F(R)$ co\"{\i}ncident avec les composantes analytiques de $F(R)$}~; on peut comparer avec~\cite{Be non rec} et~\cite{Be comp}.

	L'image par $R$ d'une composante de l'ensemble de Fatou c'est aussi une composante de l'ensemble de Fatou.
	Alors pour une composante $C$ de $F(R)$ il y a deux cas.
	Soit les composantes $R^n(C)$, pour $n \ge 0$, sont disjointes deux \`a duex~; on dit alors que $C$ est {\it errante}.
	Soit il existe $m > n \ge 0$ tel que $R^n(C) = R^m(C)$~; alors la composante $C_0 = R^n(C)$ satisfait $R^{m - n}(C_0) = C_0$ et on dit que la composante $C_0$ est {\it p\'eriodique}.

\subsection{Classification des composantes p\'eriodiques.}
Le {\it bassin d'attraction} d'un point p\'eriodique attractif $z_0$ est l'ensemble des points $x \in \pp$ satifaisant
$$
	\Delta(R^n(x), R^n(z_0)) \rightarrow 0 
		\mbox{ lorsque } n \rightarrow \infty,
$$
o\`u $\Delta$ note la distance chordale de $\pp$.
	De plus le {\it bassin d'attraction imm\'ediat} de $z_0$ est la composante analytique du bassin d'attraction qui contient $z_0$.

Le {\it domaine de quasi-p\'eriodicit\'e} de $R$ est l'interieur de l'ensemble des points de $\pp$ r\'ecurrents par $R$.

\

	Dans~\cite{these} on a caract\'eris\'e la g\'eom\'etrie des bassins d'attraction imm\'ediat et des composantes analytiques du domaine de quasi-p\'eriodicit\'e.
	En particulier les composantes analytiques du domaine de quasi-p\'eriodicit\'e sont des affino\"{\i}des ouverts (Th\'eor\`eme~3 de~\cite{these}).
\begin{thmA}[Classification des Composantes P\'eriodiques] \

	Soit $R$ une fonction rationnelle de degr\'e au moins deux et soit $C_0$ une composante p\'eriodique de $F(R)$.
	Alors il y a deux cas.
	\begin{enumerate}
		\item[1.]
			$C_0$ est un bassin d'attraction imm\'ediat.
		\item[2.]
			$C_0$ est une composante analytique du domaine de quasi-p\'eriodicit\'e.
	\end{enumerate}
\end{thmA}

	 Dans~\cite{these} {\it on a conjectur\'e que l'ensemble de Fatou est la r\'eunion (disjointe) des bassins d'attraction et des pr\'eimages du domaine de quasi-p\'eriodicit\'e}.
	D'apr\`es le Th\'eor\`eme~A cette conjecture est \'equivalente \`a la conjecture suivante.
\begin{non errance}
	Toutes les composantes de l'ensemble de Fatou sont pr\'e-p\'eriodiques.
\end{non errance}
	Dans le cas complexe un th\'eor\`eme de D.~Sullivan dit qu'il n'y a pas de composantes errantes de l'ensemble de Fatou~; voir~\cite{Su}.
	De plus il y a une caract\'erisation analogue des composantes p\'eriodiques~:
	une composante p\'eriodique est soit un bassin d'attraction (d'un point p\'eriodique attractif ou parabolique), soit un disque de Siegel ou un anneau de Herman~; voir~\cite{Su}, \cite{Mi} ou~\cite{Bea}.
	Les disques de Siegel et anneaux de Herman sont analogues au domaine de quasi-p\'eriodicit\'e, car leur r\'eunion est \'egal \`a l'interieur de l'ensemble des points r\'ecurrents.
\subsection{Une propri\'et\'e de point fixe pour les fonctions rationnelles.}
	Fixons une fonction rationnelle $R \in \C_p(z)$.

	Une partie cl\'e dans la d\'emonstration du Th\'eor\`eme~A c'est d'analyser la situation o\`u un espace analytique $X \subset \pp$ est invariant par $R$, i.e. $R(X) \subset X$.

	Dans le cas complexe on applique le Lemme de Schwarz-Pick~: $R|_X$ est, soit une isom\'etrie, soit une contraction pour la distance hyperbolique de $X$~; voir e.g.~\cite{Mi} et~\cite{Bea}.

\begin{figure}[htb]
\begin{center}
\psfig{file=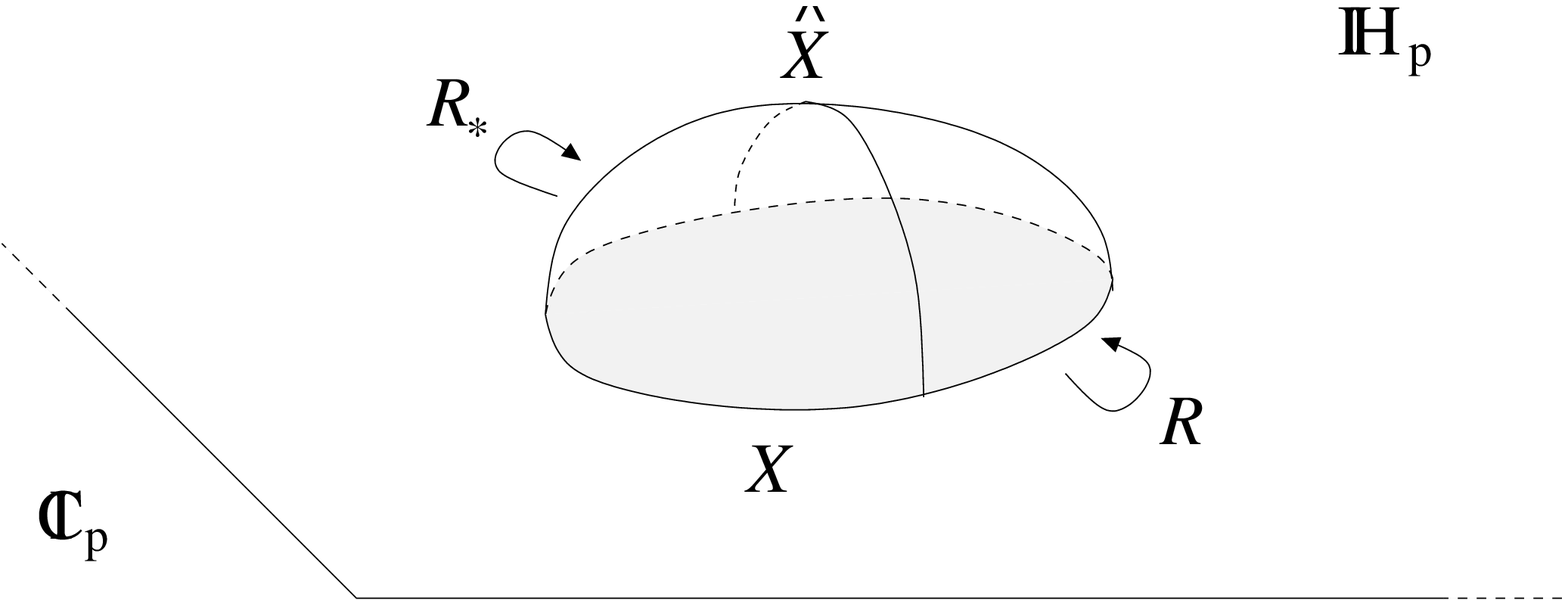, width = 3.5in}
\caption{}
\label{fig env}
\end{center}
\end{figure}

	Dans le cas $p$-adique il n'y a apparament pas d'analogue d'une distance hyperbolique d\'efinie sur $X$.
	Alors on fait un raisonement compl\`etement diff\'erent.
	{\it On consid\`ere l'action $R_*$ induite par $R$ sur l'espace hyperbolique $\H_p$}~; voir~\cite{hyp}.
	La propri\'et\'e $R(X) \subset X$ implique $R_*(\hX) \subset \hX$, o\`u $\hX \subset \H_p$ est l'{\it enveloppe convexe de $X$ dans} $\H_p$~; voir figure~\ref{fig env}.

	Alors on peut appliquer la propri\'et\'e suivante, qui est d'int\'er\^et ind\'ependent (l'ensemble $\hX$ est connexe et il contient au moins deux points).
\begin{propriete de point fixe}
	Soit $R \in \C_p(z)$ une fonction rationnelle et soit $\hX \subset \H_p$ un ensemble connexe contenant au moins deux points, tel que $R_*(\hX) \subset \hX$.
	Alors, soit $\hX$ contient un point fixe rationnel de $R_*$ ; soit il existe un point fixe attractif $z_0 \in \pp$ de $R$, tel que $\hX$ contient une demi-g\'eod\'esique issue de $z_0$.
\end{propriete de point fixe}
\noindent
(Les points rationnels de $\H_p$ sont les points de ramification de $\H_p$).

	Une cons\'equence simple de cette propri\'et\'e est la proposition suivante, laquelle est un pas important dans la preuve du Th\'eor\`eme~A.
\begin{propA}
	Soit $X \subset \pp$ un espace anlytique tel que $R(X) \subset X$.
	Alors $X$ intersecte un bassin d'attraction imm\'ediat ou le domaine de quasi-p\'eriodicit\'e de $R$.
\end{propA}
\subsection*{Remerciements.}
	Je voudrais remercier R.~Benedetto, avec qui j'ai eu des discussions r\'eli\'es a ce travail.
\section{Pr\'eliminaires.}
	Soit $p > 1$ un nombre premier, $\Q_p$ le corps des nombres $p$-adiques et soit $\C_p$ la plus petite extension compl\`ete et alg\'ebriquement close de $\Q_p$.

	On note $| \cdot |$ la norme sur $\C_p$ et $\C_p^* = \C_p - \{ 0 \}$ le groupe multiplicatif de $\C_p$.
	On appelle
\begin{eqnarray*}
	|\C_p^*| & = & \{ |z| \mid z \in \C_p^* \} \\
		& = & \{ r > 0 \mid \log_p r \mbox{ est rationnel} \}.
\end{eqnarray*}
le {\it groupe de valuation} de $\C_p^*$.
	De plus on note $\, \Delta$ \, la distance sur $\C_p$ induite par $| \cdot |$.

	On note $\co_p = \{ z \in \C_p \mid |z| \le 1 \}$ l'{\it anneau des entiers}.
	Alors $\cm_p = \{ z \in \C_p \mid |z| < 1 \}$ est un id\'eal maximal de $\co_p$.
	Le corps $\widetilde{\C}_p = \co_p / \cm_p$ est appell\'e le {\it corps r\'esiduel de} $\C_p$, qui est isomorphe \`a la fermeture algebrique $\f$ du corps fini ${\Bbb F}_p$.
	On identifie $\widetilde{\C}_p$ \`a $\f$.

	Pour $z \in \co_p$ on note $\widetilde{z}$ la projection de $z$ dans $\f$.
	Pour $\zeta \in \f$ on pose $B(\zeta) = \{ \widetilde{z} = \zeta \}$.
	Donc on a la partition,
$$
	\co_p = \sqcup_\f B(\zeta).
$$
\subsection{La droite projective.}
	On consid\`ere la droite projective $\pp$, qui est l'ensemble des droites dans $\C_p \times \C_p$ passant par $(0, 0)$.
	Pour $(x, y) \in \C_p \times \C_p - \{ (0, 0) \}$, on note $[x, y] \in \pp$ le point correspondant \`a la droite $\{ (\lambda x, \lambda y) \mid \lambda \in \C_p \}$.
	On note $\infty$ le point $[1, 0] \in \pp$ et on identifie $\pp - \{ \infty \}$ avec $\C_p$, par l'application $[\lambda, 1] \longrightarrow \lambda$.

\

	On \'etend la projection de $\C_p$ \`a $\f$ \`a une projection de $\pp = \C_p \cup \{ \infty \}$ \`a $\pf = \f \cup \{ \infty \}$, par $\widetilde{z} = \infty \in \f$, pour $z \in \{ |z| > 1 \} \cup \{ \infty \}$.
	On pose $B(\infty) = \{ \widetilde{z} = \infty \} = \{ |z| > 1 \} \cup \{ \infty \}$ et alors on a la partition
\begin{equation}\label{partition canonique}
	\pp = \sqcup_\pf B(\zeta).
\end{equation}

\

	On a une correspondance entre $PGL(2, \C_p)$ et les automorphismes de $\pp$, telle que \`a $({a \atop b} {c \atop d}) \in PGL(2, \C_p)$ correspond l'automorphsime de $\pp$ donn\'e en coordonn\'ees homog\`enes par $[x, y] \longrightarrow [ax + by, cx + dy]$.
	
	Le sous-groupe $PGL(2, \co_p)$ de $PGL(2, \C_p)$ correspond aux automorphismes qui pr\'eservent la partition $(\ref{partition canonique})$.
	De plus l'automorphisme de $\pp$ associ\'e \`a $( {a \atop c} {b \atop d}) \in PGL(2, \co_p)$ pr\'eserve {\it chaque} \'el\'ement de la partition $(\ref{partition canonique})$, si et seulement si $|a - 1| < 1$, $|d - 1| < 1$, $|b| < 1$ et $|c| < 1$.

\

	La {\it distance chordale} sur $\pp$ est d\'efinit par,
$$
	\Delta( \, [x, y], [x', y'] \, ) 
		= \frac{|xy' - x'y|}{\max \{|x|,|y|\} \max \{ |x'|,|y'| \}}
$$
o\`u en coordonn\'ees,
$$
	\Delta( z, z') = \frac{|z - z'|}{\max \{ 1, |z| \} \max \{ 1, |z'| \}},
$$
voir~\cite{Ru} ou~\cite{MS}.
	Cette distance co\"{\i}ncide avec la distance induite par $| \cdot |$ sur $\co_p$.
	De plus un automrophsime de $\pp$ pr\'eserve la distance chordale, si et seulement si correspond \`a un \'el\'ement de $PGL(2, \co_p)$.	
\subsection{Boules et couronnes.}
	Etant donn\'e $r \in |\C_p^*|$ et $a \in \C_p$ on appelle
$$
	\{ z \in \C_p \mid |z - a| < r \} \mbox{ et }
		\{ z \in \C_p \mid |z - a| \le r\}
$$
{\it boule ouverte de} $\C_p$ et {\it boule ferm\'ee de} $\C_p$, respectivement.
	Si $r \not \in |\C_p^*|$ alors ces deux ensembles co\"{\i}ncident et on l'appelle {\it boule irrationnelle de} $\C_p$.
	Notons que par d\'efinition une boule $B$ de $\C_p$ est irrationnelle si et seulement si $\diam(B) \not \in |\C_p^*|$~; en particulier {\it si $B$ est ouverte ou ferm\'ee alors $\diam(B) \in |\C_p^*|$}.

	Etant donn\'es deux boules $B$ et $B'$ de $\C_p$ il y a trois possibilit\'es~: soit $B \cap B' = \emptyset$, soit $B \subset B'$, soit $B' \subset B$.

	L'image d'une boule ouverte (resp. ferm\'ee, irrationnelle) par un automorphisme de $\C_p$ est une boule de la m\^eme nature.

\

	Une boule ouverte (resp. ferm\'ee, irrationnelle) de $\pp$ est soit une boule de $\C_p$ de la m\^eme nature, soit le compl\'ementaire d'une boule ferm\'ee (resp. ouverte, resp. irrationnelle) de $\C_p$.
	{\it Pour ce qui suit le mot boule d\'enotera une boule de $\pp$.}

	Etant donn\'es deux boules $B$ et $B'$ de $\pp$ il y a quatre possibilites~: soit $B \cap B' = \emptyset$, soit $B \subset B'$, soit $B' \subset B$, soit $B \cap B' \neq \emptyset$ et $B \cup B' = \pp$.

	Dans ce dernier cas les compl\'ementaires de $B$ et $B'$ sont disjointes~; si $B$ et $B'$ ne sont pas ferm\'ees alors on dit que $B \cap B'$ est une {\it couronne}.
	Apr\`es changement de coordonn\'ee, on peut supposer $B = \{ |z| < r \}$ et $B' = \{ |z| > r' \} \cup \{ \infty \}$ avec $r' < r$~; alors
$$
	B \cap B' = \{ z \in \C_p \mid \log_p|z| \in (\log_p r', \log_p r) \}.
$$
	On note $\mod(B \cap B') = \log_p r - \log_p r' > 0$, qui ne d\'epend pas de choix de coordonn\'ee, et on l'appelle le {\it module} de la couronne $B \cap B'$. 

	L'image d'une boule ouverte (resp. ferm\'ee, irrationnelle) par un automorphisme de $\pp$ est une boule de la m\^eme nature.
\subsection{Fonctions rationnelles.}
	Consid\'erons une fonction rationnelle $R \in \C_p(z)$ qui ne soit pas constante.
	Etant donn\'e un point $w \in \pp$ le {\it degr\'e local} de $R$ en $w$, que l'on note $\deg_R(w)$, est d\'efinit comme suit.
	On consid\`ere des coordonn\'ees tel que $w = 0$ et $R(0) = 0$.
	Alors $R$ est localement de la forme
$$
	a_dz^d + a_{d + 1} z^{d + 1} + \ldots, 
		\mbox{ o\`u } d \ge 1 \mbox{ et } a_d \neq 0 \, \, ; 
$$
on d\'efinit $\deg_R(w) = d$ et on dit que $\deg_R(w)$ est la {\it multiplicit\'e} de $w$ comme pr\'eimage de $R(w)$.
	Il n'est pas difficile de voir que $\deg_R(w)$ ne d\'epend pas de choix des coordonn\'ees.

	Pour $w \in \pp$ on a
\begin{equation}\label{preimages projective}
	\sum_{R(z) = w} \deg_R(z) = \deg(R)
\end{equation}
et pour $Q \in \C_p(z)$ on a $\deg_{Q \circ R}(w) = \deg_Q(R(w)) \cdot \deg_R(w)$.

	Etant donn\'es $X$, $Y \subset \pp$ tel que $R(X) \subset Y$ on dit que $R : X \longrightarrow Y$ est {\it de degr\'e $d$}, o\`u $d \geq 1$, si pour tout $y \in Y$
$$
        \sum_{x \in X, \, R(x) = y} \deg_R(x) = d \,\, ;
$$
de fa\c{c}on \'equivalente, tout point dans $Y$ a exactement $d$ pr\'eimages dans $X$ compt\'ees avec multiplicit\'e.

\

	Les {\it points critiques} d'une fonction rationnelle $R \in \C_p(z)$ sont les points $c \in \pp$ tel que $\deg_R(c) > 1$.
	Si $c \in \pp$ est un point critique de $R$, alors la {\it multiplicit\'e de} $c$ {\it comme point critique} de $R$ est par d\'efinition $\deg_R(c) - 1$.
	Une fonction rationnelle $R$ a $2\deg(R) - 2$ points critiques compt\'es avec multiplicit\'e. 
	C'est-\`a-dire
$$
	\sum_\pp (\deg_R(c) - 1) = 2\deg(R) - 2.
$$

	On dit qu'un point $z_0 \in \pp$ est {\it p\'eriodique} par $R$ s'il existe un entier $n \ge 1$ tel que $R^n(z_0) = z_0$.
	Alors on dit que $n$ est une p\'eriode de $z_0$ et si $n \ge 1$ est le plus petit entier avec cette propri\'et\'e on dit que $n$ est la {\it p\'eriode primitive} de $z_0$.

	Dans ce cas on appelle $\lambda = (R^n)'(z_0) \in \C_p$ le {\it multiplicateur} de $z_0$.
	Alors on dit que $z_0$ est {\it attractif}, {\it indiff\'erent} ou {\it r\'epulsif} selon que $|\lambda| < 1$, $|\lambda| = 1$ ou $|\lambda| > 1$, respectivement.
	Toute fonction rationnelle a un point fixe non r\'epulsif~; voir~\cite{Be hyp}.
\section{Bouts et espace hyperbolique}
	Dans cette section on rappelle la d\'efinition des bouts et du espace hyperbolique $\H_p$, voir~\cite{hyp} pour les r\'ef\'erences.

	L'espace hyperbolique est munit d'une distance naturelle pour laquelle est un espace m\'etrique complet~; voir Section~\ref{sec hyp}.
	On n'en aura besoin que dans la Section~\ref{propriete point fixe}, mais parfois on fera des r\'ef\'erences \`a la topologie sur $\H_p$ induite par cette distance (e.g. Lemme~\ref{detalles}).
\subsection{Bouts.}\label{bouts}
	Soit $\{ B_i \}_{i \ge 0}$ une suite croissante de boules ferm\'ees ou irrationnelles tel que $B = \cup_{i \ge 0} B_i$ soit une boule ouverte ou irrationnelle, ou soit \'egale \`a $\pp$.
	Alors $\{ B - B_i \}_{i \ge 0}$ est soit une suite decroissante de couronnes, soit une suite decroissante de boules, respectivement.
	On appelle $\{ B - B_i \}_{i \ge 0}$ {\it chaine \'evanescente}.
	Notons qu'on a
$$
	\cap_{i \ge 0} (B - B_i) = \emptyset
$$
et par cons\'equent $B - B_i \subset \C_p$, pour $i$ assez grand.
	De plus $\diam(B - B_i)$ converge vers un nombre positive lorsque $i \rightarrow \infty$.

	On dit que deux chaines \'evanescentes $\{ B - B_i \}_{i \ge 0}$ et $\{ B' - B_i' \}_{i \ge 0}$ sont {\it \'equivalentes} si pour tout $N \ge 0$ il existe $n \ge N$ tel que $B_N \subset B_n'$ et $B_N' \subset B_n$.
	Dans ce cas $B = B'$.
\begin{definition}
	Un {\bf bout} est une classe de \'equivalence de chaines \'evanescentes.
\end{definition}
	Soit $\cp$ un bout et $\{ B - B_i \}_{i \ge 0}$ une chaine \'evanescente d\'efinissante.
	Alors $B$ d\'epend seulement en $\cp$ et on note $B_\cp = B$.

	Si $B_\cp = \pp$, alors on dit que $\cp$ est un {\it bout singulier}.
	Sinon $B_\cp$ est une boule ouverte ou irrationnelle qui est d\'etermin\'ee par $\cp$.
	Si $B_\cp$ est une boule ouverte (resp. irrationnelle) alors on dit que $\cp$ est {\it rationnel} (resp. {\it irrationnel}).
	Alors $\cp \longrightarrow B_\cp$ est une correspondance entre les boules ouvertes (resp. irrationnelles) et les bouts rationnels (resp. irrationnels).

	Chaque automorphisme $\varphi$ de $\pp$ induit une bijection sur les bouts rationnels (resp. irrationnels, singuliers).
	On note cette action par $\varphi_*$.
\subsection{Points de l'espace hyperbolique.}\label{points}
	Les points de l'espace hyperbolique $\H_p$ sont des ensembles de bouts~; chaque bout appartient \`a exactement un point.
	Il y a trois types de points de $\H_p$~: les points singuliers, irrationnels et rationnels.
\subsubsection{Points singuliers.}
	Les {\it points singuliers} de $\H_p$ sont les ensembles de la forme $\cs = \{ \cp \}$, o\`u $\cp$ est un bout singulier.
\subsubsection{Points irrationnels.}
	Les {\it points irrationnels} de $\H_p$ sont les ensembles de la forme $\{ \cp, \cp' \}$ o\`u $\cp$ et $\cp'$ sont des bouts irrationnels satisfaisant $B_\cp \sqcup B_{\cp'} = \pp$.
	
	D'autre part, si $B$ est une boule irrationnelle, alors $B' = \pp - B$ est aussi une boule irrationnelle et l'ensemble ayant les bouts de $B$ et $B'$ comme \'el\'ements est un point irrationnelle de $\H_p$.
\subsubsection{Le point canonique et les points rationnels.}
	Rappelons que pour $\zeta \in \pf$ on not\'e $B(\zeta)$ la boule $\{ z \in \pp \mid \widetilde{z} = \zeta \}$~; voir Pr\'eliminaires.
	Donc on a la {\it partition canonique}
$$
	\pp = \sqcup_\pf B(\zeta).
$$
	Alors le {\it point canonique} de $\H_p$ est l'ensemble $\cs_{\can} = \{ \cp(\zeta) \}_\pf$, o\`u $\cp(\zeta)$ est le bout de $B(\zeta)$.

	De plus les {\it points rationnels} de $\H_p$ sont les ensembles de la forme $\{ \varphi_*(\cp(\zeta)) \}_\pf$, o\`u $\varphi$ est un automorphisme de $\pp$.
	En particulier $\cs_\can$ est un point rationnel.
	Notons qu'on a un param\'etrage de $\cs$ par $\pf$ qui est unique, sauf changement de coordonn\'ee projectif de $\pf$.
\subsubsection{Partitions de la droite projective}\label{partitions maximales}
	Notons que pour tout point $\cs \in \H_p$ on a la partition
\begin{equation}\label{partition point}
	\pp = \sqcup_\cs B_\cp.
\end{equation}

	Si le point $\cs$ est singulier cette partition est triviale.

	Si le point $\cs$ est rationnel (resp. irrationnel) chaque $B_\cp$, pour $\cp \in \cs$, est un boule rationnelle (resp. irrationnelle)~; on dit que $B_\cp$ est une {\it boule associ\'ee \`a} $\cs$.
	Alors on a la propri\'et\'e de maximalit\'e suivante.
	Si $B$ est une boule associ\'ee \`a $\cs$ et $D$ est une boule disjointe de $B$, alors il existe une boule $B'$ associ\'ee \`a $\cs$ qui contient $D$ (voir~\cite{hyp}).
\subsection{Espace hyperbolique.}
\begin{definition}
	L'{\bf espace hyperbolique $p$-adique}, qu'on note $\H_p$, est l'ensemble des points rationnels, irrationnels et singuliers.
	De plus on note $\H_p^\Q$ (resp. $\H_p^\R$) l'ensemble des points rationnels (resp. non singuliers) de $\H_p$.
\end{definition}
	Il est claire que chaque le groupe $PGL(2, \C_p)$ des automorphismes de $\pp$ agit sur $\H_p$, pr\'eservant $\H_p^\R$ et $\H_p^\Q$.
	Cette action est transitive sur $\H_p^\Q$ et le stabilisateur du point $\cs_\can$ correspond au groupe $PGL(2, \co_p)$.
	Par cons\'equent on a une bijection entre $PGL(2, \C_p)/PGL(2, \co_p)$ et $\H_p^\Q$.

	Etant donn\'e un automorphisme $\varphi$ de $\pp$ on note $\varphi_*$ l'action sur $\H_p$ induite par $\varphi$.
\subsubsection{Propri\'et\'e de s\'eparation.}
\begin{definition}\label{separation}
	Soit $\cs \in \H_p$ et $X \subset \pp$.
	\begin{enumerate}
		\item
	Si $\cs \in \H_p^\R$ est non singulier et si $X$ intersecte au moins deux boules associ\'ees \`a $\cs$, alors on dit que $\cs$ {\bf s\'epare} $X$ et on note $\cs \prec X$.
		\item
	Si $\cs = \{ \cp \} \in \H_p$ est singulier et si pour toute chaine \'evanescente $\{ D_i \}_{i \ge 0}$ d\'efinissant $\cp$ on a $D_i \subset X$, pour $i$ assez grand, alors on note $\cs \prec X$.
	\end{enumerate}
\end{definition}
	Soient $\cs \in \H_p$ et $X$, $Y \subset \pp$.
	Alors $\cs \prec X$ et $X \subset Y$ implique $\cs \prec Y$.
	D'autre part, si $\cs$ est singulier, $\cs \prec X$ et $\cs \prec Y$ implique $X \cap Y \neq \emptyset$.
\subsubsection{G\'eod\'esiques et demi-g\'eod\'esiques.}
	Pour $r \in \R$ on note $\cs_r$ le point associ\'e \`a la boule $\{ |z| < p^r \}$.

	Etant donn\'e deux points distincts $z, z' \in \pp$ la {\it g\'eod\'esique} de $\H_p$ joignant $z$ et $z'$, qu'on note $(z, z')$, est l'ensemble de tous les points de $\H_p^\R$ qui s\'epare $\{ z, z' \} \subset \pp$.
	Apr\`es changement de coordonn\'ee on peut supposer $z = 0$ et $z' = \infty$ et alors la g\'eod\'esique correspondant est \'egale \`a $\{ \cs_r \}_{\R}$.
	Si l'on munit $\H_p$ de sa distance naturelle alors chaque g\'eod\'esique est isom\'etrique \`a $\R$~; voir~\cite{hyp}.

	Etant donn\'e $z \in \pp$ et $\cs \in \H_p^\R$, soit $\cp \in \cs$ le bout tel que $z \in B_\cp$.
	Alors la {\it demi-g\'eod\'esique} joignant $z$ et $\cs$, qu'on note $(z, \cs)$, est l'ensemble de tous les points de $\H_p^\R$ qui s\'epare $B_\cp - \{ z \}$ et le compl\'ementaire de $B_\cp - \{ z \}$.
	Apr\`es changement de coordonn\'ee on peut supposer $z = 0$ et $\cs = \cs_{r_0}$ et alors la demi-g\'eod\'esique correspondant est \'egale \`a $\{ \cs_r \}_{r < r_0}$~; voir aussi Lemme~\ref{demi-geodesique}.
\subsubsection{Enveloppes convexes.}\label{enveloppes convexes}
	Etant donn\'e un ensemble $X \subset \pp$ on appelle
$$
	\hX = \{ \cs \in \H_p \mid \cs \prec X \}
$$
l'{\it enveloppe convexe} de $X$.
	Notons qu'il contient tous les g\'eod\'esiques joignant points distincts de $X$, mais $\hX$ peut aussi contenir des points singuliers.
\begin{lemme}\label{detalles}
	Pour une parite $X \subset \pp$ contenant au moins deux points on a les propri\'et\'es suivantes.
	\begin{enumerate}
		\item
			L'ensemble $\hX = \{ \cs \in \H_p \mid \cs \prec X \}$ contient au moins deux points et il est connexe.
		\item
			Si $\hX$ contient une demi-g\'eod\'esique issue de $z_0 \in \pp$, alors $z_0$ appartient \`a la fermeture topologique de $X$ dans $\pp$.
	\end{enumerate}
\end{lemme}
\preuve

$1.-$	L'assertion sur la connexit\'e de $\hX$ c'est le Lemme~3.8 de~\cite{hyp}.
	Comme $X$ contient au moins deux points, l'ensemble $\hX \subset \H_p$ contient une g\'eod\'esique de $\H_p$.
	Alors $\hX$ contient au moins duex points, car chaque g\'eod\'esique de $\H_p$ contient une infinit\'e de points.

\

$2.-$	Apr\`es changement de coordonn\'ee on suppose $z_0 = 0$.
	Pour $r > 0$ soit $\cs_r \in \H_p$ le point associ\'e \`a $\{ |z| < r \}$.
	Notons que si $\cs_r$ s\'epare $X$, alors $\{ |z| \le r \}$ intersecte $X$.
	Par hypoth\`ese pour $r > 0$ petit le point $\cs_r$ s\'epare $X$, donc on conclut que $0$ appartient \`a la fermeture topologique de $X$.
$\hfill \square$
\subsection{Action des fonctions rationnelles sur les bouts et les points.}
	Fixons une fonction rationnelle $R \in \C_p(z)$.
\subsubsection{Action sur les bouts.}
	La proposition suivante d\'ecrit l'action des foncitons rationnels sur les bouts, voir Proposition~4.1 de~\cite{hyp}.
\begin{proposition}\label{bout}
	Soit $\cp$ un bout rationnel (resp.  irrationnel, singulier).
	Alors il existe un bout $\cp'$ de la m\^eme nature et un entier $d \ge 1$ tel que pour toute chaine \'evanescente $\{  C_i \}_{i \ge 0}$ d\'efinissant $\cp$, il existe $N \ge 1$ tel que on ait les propri\'et\'es suivantes.
	\begin{enumerate}
		\item
			$\{ R(C_i) \}_{i \ge N}$ est une chaine \'evanescente d\'efinissant $\cp'$.
		\item
			Pour tout $i \ge N$, $R : C_i \longrightarrow R(C_i)$ est de degr\'e $d$.
	\end{enumerate}
\end{proposition}
	On note $R_*(\cp) = \cp'$ et $\deg_R(\cp) = d$.

	Le lemme suivant c'est le Lemme~4.2 de~\cite{hyp}.
\begin{lemme}\label{residu}
	Soit $\cp$ un bout non singulier.
	Alors, soit $R : B_\cp \longrightarrow B_{R_*(\cp)}$ est de degr\'e $\deg_R(\cp)$, soit $R(B_\cp) = \pp$.
	En particulier dans les deux cas on a $B_{R_*(\cp)} \subset R(B_\cp)$.
\end{lemme}
\subsubsection{Action sur les points de $\H_p$.} \

	Etant donn\'e un point singulier $\cs = \{ \cp \} \in \H_p$ on note
$$
	R_*(\cs) = \{ R_*(\cp) \} \mbox{ et } 
		\deg_R(\cs) = \deg_R(\cp) \ge 1.
$$

	Pour un point irrationnel $\cs = \{ \cp, \cp' \} \in \H_p$ on a
$$
	\deg_R(\cp) = \deg_R(\cp') \ge 1 \mbox{ et }
		\{ R_*(\cp), R_*(\cp') \} \in \H_p
$$
est un point irrationnel de $\H_p$.
	On note $\deg_R(\cs)$ et $R_*(\cs) \in \H_p$ respectivement.

	La proposition suivante d\'ecrit l'action d'une fonction rationnelle sur les points rationnels de $\H_p$~; voir~\cite{these} Proposition $2.4$.
\begin{proposition}\label{point rationnel}
	Soit $R \in \C_p(z)$ une fonction rationnelle et $\cs \in \H_p^\Q$ un point rationnel.
	Alors on a les propri\'et\'es suivantes.
          \begin{enumerate}
                \item
			Il existe un point rationnel $\cs' \in \H_p^\Q$ tel que si $\cp \in \cs$ alors $R_*(\cp) \in \cs'$.
		\item
	Consid\'erons des param\'etrages
$$
	\cs = \{ \cp(\xi) \}_{\xi \in \pf}
		\mbox{ et } \cs' = \{ \cp'(\xi) \}_{\xi \in \pf}.
$$
	Alors il existe une fonction rationnelle $\widetilde{R} \in \f(z)$ telle que pour tout $\xi \in \pf$ on a
$$
	R_*(\cp(\xi)) = \cp'(\widetilde{R}(\xi)) \mbox{ et }
		\deg_R(\cp(\xi)) = \deg_{\widetilde{R}}(\xi).
$$

	Donc pour tout $\cp' \in \cs'$
$$
        \sum_{\cp \in \cs, R_*(\cp) = \cp'} \deg_R(\cp) = \deg(\widetilde{R}).
$$
                \item
	Il existe un sous-ensemble fini $\ct \subset \cs$ tel que $R(D_\cp) = \pp$ pour tout $\cp \in \ct$ et tel que $R : D_\cp \longrightarrow D_{R_*(\cp)}$ est de degr\'e $\deg_R(\cp)$ pour tout $\cp \in \cs - \ct$~;
	dans ce dernier cas $R(D_\cp) = D_{R_*(\cp)}$.
	\end{enumerate}  
\end{proposition}
	On note $\deg_R(\cs) \ge 1$ le degr\'e de $\widetilde{R}$, qui ne d\'epend pas de choix des coordonn\'ees.
\subsubsection{Points p\'eriodiques.}
	On dit qu'un point $\cs \in \H_p$ est {\it p\'eriodique} par $R_*$ s'il existe un entier $n \ge 1$ tel que $R_*^n(\cs) = \cs$.
	Dans ce cas on dit que $\cs$ est {\it r\'epulsif} si $\deg_{R^n}(\cs) > 1$ et on dit que $\cs$ est {\it indiff\'erent} si $\deg_{R^n}(\cs)  = 1$.
	
	{\it Tout point p\'eriodique r\'epulsif est rationnel}~; voir~\cite{hyp} Proposition~5.4.
\section{Affino\"{\i}des et espaces analytiques (ouverts).}
	Un {\it affino\"{\i}de ouvert} (resp. {\it ferm\'e}) est une intersection finie non vide de boules ouvertes (resp. ferm\'ees).

	Le compl\'ementaire d'un affino\"{\i}de ouvert (resp. ferm\'e) $X \subset \pp$ est une union finie disjointe de boules ferm\'ees (resp. ouvertes) $B_1, \ldots, B_k$.
	Alors les {\it bouts de} $X$ sont par d\'efinition les bouts des boules ouvertes $\pp - B_i$.
	Donc si $\cp_1, \ldots, \cp_k$ sont les bouts d'un affino\"{\i}de ouvert $X$, on a $X = \cap B_{\cp_i}$.

	L'union de deux d'affino\"{\i}des ouvertes (resp. ferm\'ees) dont l'intersection est non vide est un affino\"{\i}de ouvert (resp. ferm\'e).
	De plus un intersection finie non vide de affino\"{\i}des ouverts (resp. ferm\'es) est un affino\"{\i}de ouvert (resp. ferm\'e).

	Un {\it espace analytique} ({\it ouvert}) est l'union (finie ou infinie) de affino\"{\i}des ferm\'es (resp. ouverts) contenant un m\^eme point.
	Comme toute affino\"{\i}de ouvert est une union croissante de affino\"{\i}des frem\'es {\it tout espace analytique ouvert est un espace analytique}.
	Mais par example un affino\"{\i}de ferm\'e est un espace analytique qui n'est pas un espace analytique ouvert.
	Les affino\"{\i}des ouverts et les couronnes sont des espaces analytiques ouverts.

	Rappelons que pour une partie ouverte $U$ de $\pp$ et $x \in U$ on a d\'efinit la {\it composante analytique de $U$ qui contient} $x$ comme l'union de tous les affino\"{\i}des ouverts contenus dans $U$ et qui contient $x$~; voir Introduction.

	Comme l'union de deux affino\"{\i}des ouverts contenant un m\^eme point est aussi un affino\"{\i}de ouvert, deux composantes analytiques qui s'intersectent co\"{\i}ncident.
	En particulier $U$ est l'union disjointe de ces composantes analytiques.
	De plus notons que toute composante analytique est un espace analytique ouvert.

\begin{remarque}
	La d\'efinition de composante analytique consid\`ere ici est plus restreint que celle introduite par R.~Benedetto dans~\cite{Be comp}, voir aussi~\cite{these}~; ici on consid\`ere affino\"{\i}des ouverts au lieu d'affino\"{\i}des ferm\'es.
	D'apr\`es le Corollaire~\ref{coincidence} (plus bas) ces deux d\'efinitions co\"{\i}ncident lorsque la composante analytique au sense de~\cite{Be comp} (i.e. avec affino\"{\i}des ferm\'es) est un espace analytique ouvert.
	Voir aussi Remarque~\ref{rem coincidence}.
\end{remarque}
\subsection{Espaces analytiques et propri\'et\'e de s\'eparation.}
	L'objectif de cette section est de montrer la propri\'et\'e suivante.
\begin{proposition}\label{espace separation}
	Soit $\cs \in \H_p$ un point non singulier et soit $X \subset \pp$ un espace analytique tel que $\cs \prec X$.
	Alors l'ensemble $\ct$ des bouts $\cp \in \cs$ tel que $B_\cp \not\subset X$ est fini.
	Si de plus $X$ est un espace analytique ouvert, alors pour tout $\cp \in \ct$ il existe une couronne $C_\cp \subset X$ ayant $\cp$ comme bout.
\end{proposition}
	La preuve de cette proposition est plus bas.
	D'abord on consid\`ere quelque corollaires.
\begin{corollaire}
	Soit $X \subset \pp$ un espace analytique ouvert.
	Alors pour tout affino\"{\i}de ferm\'e $Y \subset X$ il existe un affino\"{\i}de ouvert $Y'$ tel qu'on ait $Y \subset Y' \subset X$.
\end{corollaire}
\preuve
	Comme $Y$ est un affino\"{\i}de ferm\'e son compl\'ementaire dans $\pp$ est une r\'eunion disjointe de boules ouvertes $B_1, \ldots, B_k$.
	Soit $\cs_i$ le point (rationnel) associ\'e \`a $B_i$.
	Alors chaque boule $B_i$ est contenue dans une boule associ\'ee \`a $\cs_i$~(voir Section~\ref{partitions maximales}).
	Comme il existe une infinit\'e de boules associ\'ees \`a $\cs_i$ (car $\cs_i$ est rationnel) on a $\cs_i \prec X$.
	Par la Proposition~\ref{espace separation} pour chaque boule $B_i$ il existe une boule ferm\'ee $D_i \subset B_i$ telle que $C_i = B_i - D_i \subset X$.
	Alors $Y' = \pp - D_1 \sqcup \ldots \sqcup D_k$ est un affino\"{\i}de ouvert qui contient $Y$ et qui est contenu dans $X$.
$\hfill \square$

\

	Le corollaire suivante est une cons\'equence imm\'ediate du corollaire pr\'ec\'edent
\begin{corollaire}\label{coincidence}
	Soit $U \subset \pp$ un ouvert, soit $x \in U$ et soit $X \subset U$ l'union de tous les affino\"{\i}des {\rm ferm\'es} contenant $x$ et contenus dans $U$.
	Si l'espace analytique $X$ est un espace analytique ouvert, alors $X$ co\"{\i}ncide avec la composante analytique de $U$ qui contient $x$.
\end{corollaire}

	Le preuve de la Proposition~\ref{espace separation} d\'epend du lemme suivant.
\begin{lemme}\label{0}
	Soit $\cs \in \H_p^\R$ un point non singulier et $B' \subset \pp$ une boule telle que $\cs \prec \pp - B'$.
	Alors il existe une unique boule $B$ associ\'ee \`a $\cs$ qui intersecte $B'$ et alors $B' \subset B$.
\end{lemme}
\preuve
	Comme $\cs \prec D' = \pp - B'$ il existe des boules distincts $B_0$ et $B_1$ associ\'ees \`a $\cs$ qui intersecte $D'$.
	Apr\`es changement de coordonn\'ee on peut supposer $0 \in B_0 \cap D'$ et $\infty \in B_1 \cap D'$.
	Alors il exsite $r > 0$ tel que $B_0 = \{ |z| < r \}$ et $B_1 = \{ |z| > r \} \cup \{ \infty \}$.
	De plus $B' = \pp - D' \subset \C_p - \{ 0 \}$ et par cons\'equent $B'$ est de la forme $\{ |z - a| < s \}$ ou $\{ |z - a| \le s \}$, o\`u $0 < s \le |a|$ ou $0 < s < |a|$ respectivement.

	Si $|a| > r$ (resp. $|a| < r$, $|a| = r$), alors $B_0$ (resp. $B_1$, $\{ |z - a| < r \}$) est une boule associ\'ee \`a $\cs$ qui contient $B'$, et par cons\'equent elle est la seule boule associ\'ee \`a $\cs$ qui intersecte $B'$.
$\hfill \square$

\

{\par\noindent {\it Prueve de la Proposition~\ref{espace separation}.}
	Soit $Y \subset X$ un affino\"{\i}de ferm\'e tel que $\cs \prec Y$.
	Alors $\pp - Y$ est une r\'eunion finie disjointe de boules ouvertes $B_1, \ldots, B_n$.
	Si $\cp \in \ct$, alors $B_\cp$ intersecte un des boules $B_i$ et comme $\cs \prec Y \subset \pp - B_i$ on a $B_i \subset B_\cp$ (Lemme~\ref{0}).
	Donc l'ensemble $\ct$ est fini.

	Si l'espace analytique $X$ est ouvert, soit $Z \subset X$ un affino\"{\i}de ouvert tel que $\cs \prec Z$.
	Alors $\pp - Z$ est une r\'eunion disjointe de boules ferm\'ees $D_1, \ldots, D_m$.
	Comme dans le cas pr\'ec\'edent, pour tout bout $\cp \in \ct$ la boule $B_\cp$ contient les boules $D_i$ qui elle intersecte (Lemme~\ref{0}).
	Comme les boules $D_i$ sont ferm\'ees on peut choisir une boule ferm\'ee $D_\cp \subset B_\cp$ contenant les $D_i$ contenus dans $B_\cp$.
	Alors $C_\cp = B_\cp - D_\cp \subset Z \subset X$ est une couronne ayant $\cp$ comme bout.
$\hfill \square$
\subsection{Espaces analytiques et fonctions rationnelles.}
	L'image d'un affino\"{\i}de ouvert (resp. ferm\'e) par une fonction rationnelle est un affino\"{\i}de ouvert (resp. ferm\'e)~; voir e.g. Proposition~2.6 de~\cite{these}.
	Par cons\'equent l'image par une fonction rationnelle d'un espace analytique (ouvert) est aussi un espace analytique (ouvert).
	
	De plus on a la propri\'et\'e suivante~; voir Proposition~2.6 de~\cite{these}.
\begin{proposition}\label{preimages espaces}
	Soit $X \subset \pp$ un espace analytique (resp. affino\"{\i}de) ouvert.
	Alors $R^{-1}(X)$ a un nombre fini de composantes analytiques $X_1, \ldots, X_k$ et pour chaque $0 \le i \le k$ l'application $R : X_i \longrightarrow X$ est de degr\'e $d_i \ge 1$, o\`u les entiers $d_i$ satisfont $d_1 + \cdots + d_k = \deg(R)$.
\end{proposition}
	On aura besoin du lemme suivante dans la Section~\ref{propriete point fixe}.
\begin{lemme}\label{espaces extension}
	Soit $X \subset \pp$ un espace analytique et soit $R \in \C_p(z)$ une fonction rationnelle.
	Alors si $\cs \in \H_p$ satifait $\cs \prec X$, on a $R_*(\cs) \prec R(X)$.
\end{lemme}
\preuve
	Soit $\cs \in \H_p$ tel que $\cs \prec X$.

	Si $\cs = \{ \cp \}$ est singulier, soit $\{ D_i \}_{i \ge 0}$ une chaine \'evanescente d\'effinisant $\cp$.
	On suppose que $\{ R(D_i) \}_{i \ge 0}$ est une chaine \'evanescente d\'effinisant $R_*(\cp)$ (voir Proposition~\ref{bout}).
	Par hypothese on a $D_i \subset X$ pour $i \ge 0$ assez grand et par cons\'equent $R(D_i) \subset R(X)$.

	Supposons $\cs$ non singulier.
	Soit $\ct$ l'ensemble fini des bouts $\cp \in \cs$ tel que $B_\cp \not \subset X$ (Proposition~\ref{espace separation}).
	D'apr\`es le Lemme~\ref{residu}, pour tout bout $\cp \in R_*(\cs)$ qui n'appartient pas \`a l'ensemble fini $R_*(\ct)$, on a $B_\cp \subset R(X)$.
$\hfill \square$
\section{Points p\'eriodiques r\'epulsifs.}\label{sec repulsif}
	Fixons une fonction rationnelle $R \in \C_p(z)$ de degr\'e au moins duex.

	D'apr\`es un th\'eor\`eme de L.C.~Hsia (\cite{Hs}), pour tout ouvert $X \subset \pp$ contenant un point p\'eriodique r\'epulsif de $R$, on a
\footnote{En effet, ceci est vraie pour tout ouvert $X$ qui intersecte l'ensemble de Julia, voir Corollary~2.4 et Remark~2.7 de~\cite{Hs}.}
\begin{equation*}
	\pp - E(R) \subset \cup_{n \ge 0} R^n(X),
\end{equation*}
o\`u l'ensemble $E(R) \subset \pp$ cointient au plus deux points (il est appel\'e l'{\it ensemble exceptionnel} de $R$, voir Section~\ref{exceptionnels}).

	La proposition suivante donne un analogue de cette propri\'et\'e.
\begin{proposition}\label{point repulsif}
	Soit $X \subset \pp$ un espace analytique ouvert.
	S'il existe un point p\'eriodique r\'epulsif $\cs \in \H_p$ (pour l'action de $R$) tel que $\cs \prec X$, alors on a
$$
	\pp - E(R) \subset \cup_{n \ge 0} R^n(X).
$$
\end{proposition}
	La preuve de cette proposition est dans la Section~\ref{repulsif}.
	Dans la Section~\ref{exceptionnels} on consid\`ere des propri\'et\'es de l'ensemble exceptionnel.
\subsection{Ensembles exceptionnels.}\label{exceptionnels}
\begin{definition}\label{defexceptionnel}
	Soit $R \in \C_p(z)$ (resp. $\tr \in \f(z)$) une fonciton rationnelle de degr\'e au moins deux.
	Alors on dit qu'un point de $\pp$ (resp. $\pf$) est {\bf exceptionnel}, s'il a un nombre fini de pr\'eimages par les it\'er\'es de $R$ (resp. $\tr$).
	On note $E(R) \subset \pp$ l'ensemble des points exceptionnels de $R$.
\end{definition}
\begin{proposition}\label{propexceptionnel}
	Soit $R \in \C_p(z)$ (resp. $\tr \in \f(z)$) une fonction rationnelle de degr\'e au moins deux.
\begin{enumerate}
	\item
		Tout point exceptionnel $x$ est p\'eriodique et on a
$$
	\deg_R(x) = \deg(R) > 1 \mbox{ (resp. $\deg_{\tr}(x) = \deg(\tr) > 1$)}.
$$
	\item
		$R$ a au plus $2$ points exceptionnels.
\end{enumerate}
\end{proposition}
\preuve

$1.-$	Soit $E$ l'ensemble fini de toutes les pr\'eimages de $x$.
	Alors $R$ (resp. $\tr$) induit une bijection sur $E$.
	Par cons\'equent tout point dans $E$ a une unique pr\'eimage par $R$ (resp. $\tr$) dans $\pp$ (resp. $\pf$), d'o\`u l'assertion 1.

\

$2.-$	Par $1$, tout point exceptionnel est un point critique de multiplicit\'e $\deg(R) - 1$.
	Comme $R$ a $2 \deg(R) - 2$ points critiques compt\'es avec multiplicit\'e, on conclut que $R$ a au plus $2$ points exceptionnels.
$\hfill \square$
\subsection{Preuve de la Proposition~\ref{point repulsif}.}\label{repulsif}
	Le lemme suivant d\'epend du Lemme~\ref{couronne1} dans l'Appendice~1.
\begin{lemme}\label{bout repulsif}
	Soit $\cp$ un bout tel que $R_*(\cp) = \cp$ et $\deg_R(\cp) > 1$.
	Pour une couronne $C \subset \pp$ ayant $\cp$ comme bout, il y a duex cas.
	\begin{enumerate}
		\item
			Il existe $n \ge 1$ tel que $B_\cp \subset R^n(C)$.
		\item
			Il existe $x \in B_\cp$ tel que $B_\cp - \{ x \} \subset \cup_{n \ge 0} R^n(C)$.
	\end{enumerate}
\end{lemme}
\preuve
	Supposons que pour tout $n \ge 1$ la boule $B_\cp$ n'est pas contenue dans $R^n(C)$.	

	Pour $n \ge 0$ on d\'efinit inductivement une couronne $C_n \subset R^n(C)$ ayant $\cp$ comme bout.
	Posons $C_0 = C$ et supposons que $C_n$ est d\'ej\`a d\'efinit.
	Comme $C_n \subset R^n(C)$ et $B_\cp \not \subset R^{n + 1}(C)$, on a $B_\cp \not \subset R(C_n)$.
	Donc par le Lemme~\ref{couronne1}
$$
	C_{n + 1} = R(C_n) \cap B_\cp \subset R^{n + 1}(C)
$$ 
est une couronne ayant $\cp$ comme bout.

	De plus on a $\mod(C_{n + 1}) \ge \mod(C_n) + m_0$, o\`u $m_0 = \min\{ m, \mod(C) \}$ et $m > 0$ est donn\'e par le Lemme~\ref{couronne1}.
	Donc $\mod(C_n) \ge n \cdot m_0 \rightarrow \infty$ lorsque $n \rightarrow \infty$.
	Par cons\'equent le diam\`etre chordale de la boule $B_n = B_\cp - C_n$ t\'end vers z\'ero lorsque $n \rightarrow \infty$.
	
	Quite \`a r\'eduire $C$ on suppose $C \subset R(C)$.
	Alors $C_n \subset C_{n + 1}$ et par cons\'equent $\{ B_n \}_{n \ge 0}$ est une suite decroissante.
	Donc $\cap_{n \ge 0} B_n$ consiste d'un point (car $\C_p$ est complet) qu'on note $x$.
	On a
$$
	B_\cp - \{ x \} = \cup_{n \ge 0} C_n \subset \cup_{n \ge 0} R^n(C).
\hfill \square
$$

\

{\par\noindent {\it Preuve de la Proposition~\ref{point repulsif}.}}
	Quitte \`a remplacer $R$ par un it\'er\'e on suppose que $\cs$ est fix\'e par l'action de $R$.
	D'apr\`es le Proposition~\ref{espace separation} il existe un ensemble fini $\ct \subset \cs$ tel que $B_\cp \subset X$ pour $\cp \in \cs - \ct$ et tel que pour $\cp \in \ct$ il existe une couronne $C_\cp \subset X$ ayant $\cp$ comme bout.

	Soit $\cp \in \ct$.
	S'il existe $\cp' \in \cs - \ct$ et $n \ge 1$ tel que $R^n_*(\cp') = \cp$, on a $B_\cp \subset B_{R_*^n(\cp')} \subset R^n(X)$ (Lemme~\ref{residu}).
	Sinon tous les pr\'eimages de $\cp$ dans $\cs$ par les it\'er\'es de $R_*$ sont dans $\ct$.
	Comme l'ensemble $\ct$ est fini, le bout $\cp$ correspond \`a un point exceptionnel de $R_*|_\cs$ (Proposition~\ref{point rationnel} et D\'efinition~\ref{defexceptionnel}).
	Par cons\'equent $\cp$ est p\'eriodique et on a $\deg_R(\cp) > 1$ (Proposition~\ref{propexceptionnel}).
	Par le Lemme~\ref{bout repulsif}, soit il existe $n \ge 1$ tel que $B_\cp \subset R^n(B_\cp)$, soit il existe $x_\cp \in B_\cp$ tel que
$$
	B_\cp - \{ x_\cp \} \subset \cup_{n \ge 0}R^n(C_\cp) 
		\subset \cup_{n \ge 0} R^n(X).
$$
	Dans les deux cas $\pp - \cup \{ x_\cp \} \subset \cup_{n \ge 0} R^n(X)$.
	Si $x_\cp$ n'est pas un point exceptionnel alors il existe $m_\cp \ge 1$ tel que $x_\cp \in R^{m_\cp}(\pp - \cup \{ x_\cp \}) \subset \cup_{n \ge 0} R^n(X)$.
	On conclut $\pp - E(R) \subset \cup_{n \ge 0} R^n(X)$.
$\hfill \square$
\section{Bassins d'attraction et le domaine de quasi-p\'eriodicit\'e.}
	Fixons une fonction rationnelle $R \in \C_p(z)$ de degr\'e au moins deux.
	Rappelons les d\'efinitions suivantes, voir Introduction et~\cite{these}.

\

\noindent
- Etant donn\'e un point p\'eriodique attractif $z_0 \in \pp$ de $R$, son {\it bassin d'attraction} c'est l'ensemble des points $x \in \pp$ satifaisant
$$
	\Delta(R^n(x), R^n(z_0)) \rightarrow 0 
		\mbox{ lorsque } n \rightarrow \infty.
$$
	Alors le {\it bassin d'attraction imm\'ediat} de $z_0$ est la composante analytique du bassin d'attraction qui contient $z_0$.

\

\noindent
- Le {\it domaine de quasi-p\'eriodicit\'e} ${\cal E}(R) \subset \pp$ de $R$ est l'interieur de l'ensemble des points de $\pp$ r\'ecurrents par $R$.

\

	L'objectif de cette section est de montrer la proposition suivante.
	Rappelons que $E(R) \subset \pp$ note l'ensemble exceptionnel de $R$, qui contient au plus deux points~; voir Section~\ref{exceptionnels}.
\begin{proposition}\label{maximalite0}
	Soit $C$ un bassin d'attraction imm\'ediat ou une composante analytique du domaine de quasi-p\'eriodicit\'e.
	Si $X \subset \pp$ est un espace analytique ouvert qui intersecte $C$ et qui n'est pas contenu dans $X$, on a
$$
	\pp - E(R) \subset \cup_{n \ge 0} R^n(X).
$$
\end{proposition}
	La preuve de cette proposition est bas\'ee sur les Th\'eor\`emes~2 et ~3 de~\cite{these} et sur la Proposition~\ref{point repulsif} dans la Section~\ref{sec repulsif}.
	Le cas des composantes analytiques du domain de quasi-p\'eriodicit\'e est dans la Section~\ref{qp} et celui des bassin d'attraction imm\'eidats dans la Section~\ref{bassins}.

	D'abord on montre la proposition suivante.
\begin{proposition}\label{maximalite1}
	Soit $C$ un bassin d'attraction imm\'ediat ou une composante analytique du domaine de quasi-p\'eriodicit\'e.
	Alors $\cup_{n \ge 0} R^n(C)$ omet au moins $3$ points de $\pp$.
\end{proposition}
	Cette proposition est une cons\'equence imm\'ediate du lemme suivant.
\begin{lemme}\label{maximalite2} \

\begin{enumerate}
	\item
		Le compl\'ementaire d'un bassin d'attraction d'un cycle est infini.
	\item
		Pour tout $n \ge 0$ le compl\'ementaire de $R^{-n}(\ce(R))$ est infini.
\end{enumerate}
\end{lemme}
\preuve

$1.-$	Toute fonction rationnelle a une infinit\'e de points p\'eriodiques (Remarque~\ref{Baker}) et un bassin d'attraction d'un cycle contient un nombre fini de points p\'eiodiques (les \'el\'ements du cycle).

\

$2.-$	On suppose $\ce(R) \neq \emptyset$~; alors $\ce(R)$ est infini car il est un ouvert. 
	Pour $n \ge 0$ posons $\ce_n(R) = R^{-n}(\ce(R))$.
	Comme $R$ est injective sur $\ce_0(R) = \ce(R)$ (car par d\'efinition tout point de $\ce(R)$ est r\'ecurrent par $R$), l'ensemble $\ce_1(R) - \ce_0(R)$ est infini.
	Donc pour tout $n \ge 0$ l'ensemble $\ce_{n + 1}(R) - \ce_n(R) = R^{-n}(\ce_1(R) - \ce_0(R))$ est infini.
$\hfill \square$
\subsection{Domaine de quasi-p\'eriodicit\'e.}\label{qp}
	Rappelons qu'on note $R_*$ l'action sur $\H_p$ induite par $R$.
	Le th\'eor\`eme suivant est une reformulation du Th\'eor\`eme~3 de~\cite{these}, voir Remarque~\ref{rem coincidence}.
\begin{thm }[\cite{these} Th\'eor\`eme~3]
	Soit $R \in \C_p(z)$ une fonction rationnelle de degr\'e au moins deux.
	Alors chaque composante analytique $C \subset \pp$ du domaine de quasi-p\'eriodicit\'e de $R$ est un affino\"{\i}de ouvert et chaque point de $\H_p$ contenant un bout de $C$ est un point p\'eriodique r\'epulsif de $R_*$.
\end{thm }
{\par\noindent {\it Preuve de la Proposition~\ref{maximalite0} dans le cas des composantes analytiques du domaine de quasi-p\'eriodicit\'e.}}

	Soient $\cp_1, \ldots, \cp_k$ les bouts de $C$ et $\cs_1, \ldots, \cs_k$ les points correspondants.
	Comme $X$ n'est pas contenu dans $C = \cap B_{\cp_i}$ il existe $1 \le i \le k$ tel que $X$ intersecte le compl\'ementaire de $B_{\cp_i}$.
	D'autre part $X$ intersecte $C \subset B_{\cp_i}$, donc $\cs_i \prec X$.
	Alors l'asseriton suit de la Proposition~\ref{point repulsif}.
$\hfill \square$
\subsection{Bassins d'attraction imm\'ediat.}\label{bassins}
	Soit $z_0 \in \pp$ un point p\'eriodique attractif de $R$~; quitte \`a changer $R$ par un it\'er\'e on suppose que $z_0$ est fix\'e par $R$.
\begin{thm }[\cite{these} Th\'eor\`erme~2]
	Soit $C \subset \pp$ le bassin d'attraction imm\'ediat de $z_0$ et soit $D$ une boule ouverte telle que $z_0 \in R(D) \subset D$.
	Pour $n \ge 0$ notons $X_n$ la composante analytique de $R^{-n}(D)$ qui contient $z_0$.
	Alors $C = \cup_{n \ge 0} X_n$.
	De plus il y a deux cas
	\begin{enumerate}
		\item
			$C$ est une boule ouverte.
		\item
			$C$ est un espace analytique de type Cantor.
	\end{enumerate}
\end{thm }
	On aura pas besion de d\'efinir les espaces analytiques de type Cantor, on renvoie le lecteur \`a~\cite{these}.

	Notons que chaque $X_n$ est un affino\"{\i}de ouvert (Proposition~\ref{preimages espaces}) et par cons\'equent $C = \cup_{n \ge 0} X_n$ est un espace analytique ouvert~; voir Remarque~\ref{rem coincidence}.

\

{\par\noindent {\it Preuve de la Proposition~\ref{maximalite0} dans le cas des bassins d'attraction.}

$1.-$	Supposons que $X_n$ est une boule pour tout $n \ge 0$~: ce cas correspond au cas o\`u $C$ est une boule ouverte.

	Comme $R : C \longrightarrow C$ est de degr\'e au moins deux (cf. Corollaire~4.16 de~\cite{these}) le point de $\H_p$ contenat le bout de $C$ est un point fixe r\'epulsif pour l'action de $R$.
	Alors la d\'emonstration est simlaire au cas des composantes analytiques du domaine de quasi-p\'eriodicit\'e.

\

$2.-$	Supposons alors qu'il existe $n_0 \ge 1$ tel que $X_n$ n'est pas une boule~: ce cas correspond au cas o\`u $C$ est de ``type Cantor''.

	Soit $X \subset \pp$ un espace analytique ouvert qui intersecte $C$ et le compl\'ementaire de $C$~; voir Remarque~\ref{bassin non trivial}.
	On se ram\`ene au cas o\`u $X$ est un affino\"{\i}de ouvert.
	On montrera $\pp - \{ z_0 \} \subset \cup_{n \ge 0} R^n(X)$, ce qui implique l'assertion desir\'ee.

	Etant donn\'e un bout non singulier $\cp$ on note $\tB_\cp = \pp - B_\cp$, qui est une boule ferm\'ee ou irrationnelle.

\

	Sans perte de g\'en\'eralit\'e on suppose $n_0 = 1$, donc $X_1$ a au moins deux bouts.
	Rappelons que $X_0 = D$ est une boule ouverte, soit $\cp_0$ son bout.
	Comme $X_0 \subset X_1$, pour tout bout $\cq$ de $X_1$ on a $\tB_\cq \subset \tB_{\cp_0}$.
	Comme pour tout bout $\cp$ de $X_n$ on a $R_*^n(\cp) = \cp_0$, il exsite au moins deux bouts $\cq$, $\cq'$ de $X_{n + 1}$ tel que $\tB_{\cq}, \tB_{\cq'} \subset \tB_\cp$.

	Comme $X \not \subset C$, pour chaque $n \ge 0$ l'affino\"{\i}de $X$ intersecte l'ensemble
$$
	\pp - X_n = \sqcup_{\mbox{ bouts de } X_n} \tB_\cq.
$$
	Donc on peut trouver inductivement une suite $\{ \cp_n \}_{n \ge 0}$ de bouts tel que pour $n \ge 0$~: $\cp_n$ soit un bout de $X_n$, $\tB_{\cp_{n + 1}} \subset \tB_{\cp_n}$ et telle que $\tB_{\cp_n}$ intersecte $X$, pour $n \ge 0$.

	Comme $\cap_{n \ge 0} \tB_{\cp_n}$ est disjoint de $C$ et $X$ intersecte chaque $\tB_{\cp_n}$, il existe $N \ge 0$ tel que pour $n \ge N$, $\tB_{\cp_n} - \tB_{\cp_{n + 1}}$ intersecte $X$.
	Donc pour $n > N$ toute composante $D_0$ de $\pp - X$ est soit disjointe de $\tB_{\cp_n} - \tB_{\cp_{n + 1}}$, soit contenue dans $\tB_{\cp_n} - \tB_{\cp_{n + 1}}$.
	Par cons\'equent il existe $n > N$ tel que $\tB_{\cp_n} - \tB_{\cp_{n + 1}} \subset X$.
	
	Soit $\cq$ un bout de $X_{n + 1}$ tel que $\tB_\cq \subset \tB_{\cp_n}$ et $\cq \neq \cp_{n + 1}$.
	On a $\tB_\cq \subset \tB_{\cp_n} - \tB_{\cp_{n + 1}} \subset X$.
	Comme $\cq$ est un bout de $X_{n + 1}$ on a $R_*^{n + 1}(\cq) = \cp_0$ et par cons\'equent
$$
	\pp - D = \tB_{\cp_0} 
		\subset R^{n + 1}(\tB_{\cq}) \subset R^{n + 1}(X)~;
$$
voir Lemme~\ref{residu}.
	Donc pour $k \ge 0$ on a $\pp - R^k(D) \subset R^{k + n + 1}(X)$ (cf. Lemme~\ref{residu}).
	Quitte \`a r\'eduire $D$ on peut supposer $\cap_{n \ge 0} R^k(D) = \{ z_0 \}$, donc $\pp - \{ z_0 \} \subset \cup_{m \ge 0} R^m(X)$.
$\hfill \square$
\begin{remarque}\label{rem coincidence}
	Dans la d\'efinition de composante analytique dans~\cite{these} (comme dans~\cite{Be comp}), on consid\`ere des affino\"{\i}des ferm\'es au lieu des affino\"{\i}des ouverts.
	Dans les Th\'eor\`emes~2 et~3 de~\cite{these} les composantes analytiques (au sense de~\cite{these}) en question sont des espaces analytiques ouverts or le Corollaire~\ref{coincidence} montre que ces espaces analytiques sont aussi des composantes analytiques au sense qu'on consid\`ere ici.
\end{remarque}
\begin{remarque}\label{bassin non trivial}
	Contrairement aux cas pr\'ec\'edents, dans ce cas on n'utilise pas que $X$ soit un espace analytique {\rm ouvert}~; la preuve s'applique sans changement si $X$ est un espace analytique quelconque.
	Donc on obtient la propri\'et\'e plus forte suivante~: si $C \subset \pp$ est un bassin d'attraction imm\'ediat de Type Cantor de $R$ et si $X \subset \pp$ est un espace analytique intersectant $C$ et son compl\'ementaire, alors on a $\pp - E(R) \subset \cup_{n \ge 0} R^n(X)$.
\end{remarque} 
\section{L'ensemble de Fatou et ces composantes.}\label{composantes fatou}
	Fixons une fonction rationnelle $R \in \C_p(z)$ de degr\'e au moins deux.
	Rappelons qu'on a d\'efinit l'{\it ensemble de Fatou} $F(R) \subset \pp$ de $R$ comme l'ensemble de tous les points ayant un voisinage $X \subset \pp$ satisfaisant
\begin{equation}\label{ev sur}
	\cup_{n \ge 0} R^n(X) \mbox{ omet au moins $3$ points de $\pp$.}
\end{equation}
	Notons que $F(R)$ est ouvert par d\'efinition.
	L'ensemble $J(R) = \pp - F(R)$ est appel\'e l'{\it ensemble de Julia} de $R$.

	La Proposition~\ref{definitions} (plus bas) montre que cette d\'efinition de l'ensemble de Fatou co\"{\i}ncide avec la d\'efinition introduite par L.C.~Hsia dans~\cite{Hs}.
	En particulier on a $R^{-1}(F(R)) = F(R)$, $R^{-1}(J(R)) = J(R)$ et pour tout entier $n \ge 1$ on a $F(R^n) = F(R)$ et $J(R^n) = J(R)$~; voir~\cite{Hs} ou~\cite{these}.

\

	Rappelons que la {\it composante de l'ensemble de Fatou qui contient} $x \in F(R)$ est l'union de tous les affino\"{\i}des ouverts $X \subset F(R)$ contenant $x$ et qui satisfont ($\ref{ev sur}$).
	Notons que toute composante de l'ensemble de Fatou est un espace analytique ouvert.

	La proposition suivante implique que duex composantes de l'ensemble de Fatou ayant une intersection non vide sont egaux.
	En particulier $F(R)$ est l'union disjointe de ces composantes.
\begin{proposition}\label{partition fatou}
	Soient $X, X' \subset \pp$ des affino\"{\i}des ouverts satisfaisant $(\ref{ev sur})$.
	Si $X \cap X' \neq \emptyset$, alors $X \cup X'$ aussi satifait $(\ref{ev sur})$.
\end{proposition}
\preuve
	Soit $x \in X \cap X' \subset F(R)$.
	Alors il y a trois cas.
	\begin{enumerate}
		\item
			$x$ appartient \`a un bassin d'attraction.
			La Proposition~\ref{maximalite0} implique que $X \cup X'$ est contenu dans le bassin du cycle correspondant.
			Par cons\'equent $\cup_{n \ge 0} R^n(X \cup X')$ est contenu dans le m\^eme bassin.
			Alors l'assertion suit de la partie 1 du Lemme~\ref{maximalite2}.
		\item
			Il existe $k \ge 0$ tel que $R^k(x) \in \ce(R)$.
			La Proposition~\ref{maximalite0} implique $X \cup X' \subset R^{-k}(\ce(R))$ et donc $\cup_{n \ge 0} R^n(X \cup X') \subset R^{-k}(\ce(R))$.
			Alors l'assertion suit de la partie 2 du Lemme~\ref{maximalite2}.
		\item
			$x$ n'appartient pas \`a un bassin d'attraction et $R^n(x) \not \in \ce(R)$ pour tout $n \ge 0$.
		La Proposition~\ref{maximalite0} implique alors que $\cup_{n \ge 0} R^n(X \cup X')$ est disjoint de $\ce(R)$.
		Si $\ce(R) \neq \emptyset$, alors $\ce(R)$ est infini, car il est ouvert.
		Donc l'assertion suit dans ce cas.

		Supposons $\ce(R) = \emptyset$.
		Comme toute fonction rationnelle a au moins un point fixe non r\'epulsif (cf.~\cite{Be hyp}) il existe un point fixe attractif $z_0 \in \pp$ de $R$.
		Soit $\ca$ le bassin de $z_0$, qui est infini car il est ouvert.
		Comme par hypoth\`ese $x \not \in \ca$ la Proposition~\ref{maximalite0} implique que $X \cup X'$ est disjoint de $\ca$.
		Comme $R^{-1}(\ca) = \ca$ on conclut que $\cup_{n \ge 0} R^n(X \cup X')$ est disjoint de $\ca$.
	\end{enumerate}
$\hfill \square$
\begin{proposition}\label{definitions}
	Pour un point $x \in \pp$ les propri\'et\'es suivantes sont \'equivalentes.
	\begin{enumerate}
		\item
			Il existe un voisinage $X \subset \pp$ de $x$ tel que $\cup_{n \ge 1} R^n(X)$ omet au moins $3$ points de $\pp$.
		\item
			Il existe un voisinage $X \subset \pp$ de $x$ tel que la famille $\{ R^n|_X \}_{n \ge 1}$ soit uniform\'ement lipschitzienne pour la distance chordale.
	\end{enumerate}
\end{proposition}
\preuve
	L'implication $1 \Rightarrow 2$ (avec le m\^eme voisinage $X$) suit du Main Theorem de~\cite{Hs}.

	Supposons alors que la propri\'et\'e $2$ est satisfaite et soit $\delta \in (0, 1)$.
	Alors on peut choisir une boule $X$ contenant $x$, assez petite telle que le diam\`etre chordale de $X_n = R^n(X)$, pour $n \ge 0$, soit strictement plus petit que $\delta < 1 = \diam(\pp)$.
	Alors $X_n \neq \pp$ et par cons\'equent $X_n$ est une boule (Lemme~\ref{residu}).
	De plus $X_i \cap X_j \neq \emptyset$ implique $X_i \subset X_j$ ou $X_j \subset X_i$.
	Alors il y a deux cas.
	\begin{enumerate}
		\item
			Les boules $X_n$, pour $n > 0$, sont disjointes de $X_0$.
			Alors $\cup_{n \ge 0} X_n$ est disjoint de $R^{-1}(X_0)$, qui contient au moins $3$ points car $X_0$ est ouvert.
		\item
			Il existe $k > 0$ tel que $X_k$ intersecte $X_0$.
			Si $X_k \subset X_0$ alors $\cup_{n \ge 0} X_n = X_0 \cup \ldots \cup X_{k - 1}$ et une r\'eunion finie de boules de diam\`etre chordale plus petit que 1~; par cons\'equent son compl\'ementiare contient au moins $3$ points.

			D'ature part, si $X_0 \subset X_k$, alors pour $n \ge 0$ on a $X_n \subset X_{n + k}$.
			Donc $D_n = \cup_{j \ge 0} X_{n + jk}$ est une boule de diam\`etre chordale au plus \'egale \`a $\delta < 1$.
			Par cons\'equent le compl\'ementaire de $\cup_{n \ge 0} X_n = D_0 \cup \ldots \cup D_{k - 1}$ contient au moins $3$ points.
	\end{enumerate}
$\hfill \square$
\subsection{Propri\'et\'es des composantes de l'ensemble de Fatou.}
	Fixons une fonction rationnelle $R \in \C_p(z)$ de degr\'e au moins deux.
	
	Notons que tout espace analytique ouvert $C \subset \pp$ tel que $\cup_{n \ge 0} R^n(C)$ omet au moins $3$ de $\pp$, est contenu dans une composante de l'ensemble de Fatou. 
\begin{proposition}
	Si $F(R) \neq \pp$, alors les composantes de $F(R)$ co\"{\i}ncident avec les composantes analytiques de $F(R)$.
\end{proposition}
\preuve
	Par d\'efinition toute composante de $F(R)$ est contenue dans une composante analytique de $F(R)$.
	D'autre part, par hypoth\`ese $J(R) = \pp - F(R)$ est non vide et donc il est parfait, voir Theorem~2.9 de~\cite{Hs}.
	En particulier $J(R)$ est infini.
	Comme $R^{-1}(J(R)) = J(R)$, pour tout $X \subset F(R)$ l'ensemble
$$
	J(R) \subset \pp - \cup_{n \ge 0} R^n(X)
$$
contient au moins $3$ points.
$\hfill \square$
\begin{proposition}\label{invariance}
	Soit $C$ une composante de $F(R)$.
	Alors chaque composante analytique de $R^{-1}(C)$ est une composante de $F(R)$.
\end{proposition}
	Le corollaire suivant est imm\'ediat.
\begin{corollaire}
	L'image d'une composante de $F(R)$ par $R$ est aussi une composante de $F(R)$.
\end{corollaire}
	On montre aussi la propri\'et\'e suivante.
\begin{proposition}\label{Fatou itere}
	Pour $n \ge 1$ les composantes de $F(R^n)$ co\"{\i}ncident avec les composantes de $F(R)$.
\end{proposition}
	La preuve des Proposition~\ref{invariance} et~\ref{Fatou itere} sont plus bas~; on consid\`ere d'abord des lemmes.
\begin{lemme}\label{espace errant}
	Soit $C \subset \pp$ un ensemble contenant au moins trois poins.
	Si pour tout $n \ge 1$, $R^n(C)$ est disjoint de $C$, alors $\cup_{n \ge 0} R^n(C)$ omet au moins $3$ points de $\pp$.
\end{lemme}
\preuve
	L'ensemble $\cup_{n \ge 0} R^n(C)$ est disjoint de $R^{-1}(C)$, qui contient au moins $3$ points.
$\hfill \square$

\

	Notons que par d\'efintion, pour toute composante $C \subset \pp$ de l'ensemble de Fatou, $\pp - C$ contient au moins $3$ points.
\begin{lemme}\label{espace periodique}
	Soit $C \subset \pp$ un ensemble tel que $\pp - C$ contient au moins $3$ points.
	Supposons qu'il existe $m \ge 1$ tel que $R^m(C) \subset C$ et tel que les ensembles $C, R(C), \ldots, R^{m - 1}(C)$ soient disjoints deux \`a deux.
	Alors $\cup_{n \ge 0} R^n(C)$ omet au moins $3$ points de $\pp$.
\end{lemme}
\preuve
	On suppose que $m \ge 1$ est le plus petit entier avec cette propri\'et\'e.
	Si $m = 1$ alors le lemme est trivial, donc on suppose $m > 1$.
	Notons que tout point p\'eriodique de $R$ dans $C'$ a une p\'eriode primitive multiple de $m > 1$.

	Comme pour tout $k > 3$ il existe un point p\'eriodique de p\'eriode primitive $k$ de $R$ (voir Remarque~\ref{Baker} plus bas), on conclut que le compl\'ementaire de $\cup_{n \ge 0} R^n(C) = C \sqcup \ldots \sqcup R^{m - 1}(C)$ dans $\pp$ continent une infinit\'e de points.
$\hfill \square$

\

{\par\noindent {\it Preuve de la Proposition~\ref{invariance}.}}
	Notons que par d\'efinition de $F(R)$ l'image d'une composante de $F(R)$ est contenue dans une composante de $F(R)$.

	Soit $C'$ une composante analytique de $R^{-1}(C)$, il suffit de montrer que $\cup_{n \ge 0} R^n(C')$ omet au moins trois points de $\pp$.
	Pour $n \ge 0$ soit $C_n$ la composante de $F(R)$ contenant $R^n(C)$.
	Notons que toute composante de $F(R)$ qui intersecte $C'$ est contenue dans $C'$.
	
	Alors soit pour tout $n \ge 0$ la composante $C_n$ est disjointe de $C'$ et alors l'assertion suit de la Proposition~\ref{espace errant}~; soit il existe $k \ge 0$ tel que $C_k \subset C'$ et alors la proposition suit de la Proposition~\ref{espace periodique}.
$\hfill \square$

\

{\par\noindent {\it Preuve de la Proposition~\ref{Fatou itere}.}}
	Rappelons que $F(R^n) = F(R)$.
	Notons que toute composante de $F(R)$ est contenue dans une composante de $F(R^n)$.
	D'autre part l'image par $R$ d'une composante de $F(R^n)$ est contenue dans une composante de $F(R^n)$ (cf. Proposition~\ref{invariance}).
	Alors l'assertion suit des Lemmes~\ref{espace errant} et~\ref{espace periodique}.
$\hfill \square$
\begin{remarque}\label{Baker}
	Un th\'eor\`eme de I.N.~Baker dit que toute fonction rationnelle de degr\'e au moins duex a des points p\'eriodiques de tous les p\'eriodes (primitives) plus grands que 3~; voir Theorem~6.2.2 de~\cite{Bea}, o\`u l'\'ennonc\'e est pour les fonctions rationnelles \`a coefficientes dans $\C$, mais la d\'emonstration s'applique au cas de $\C_p$ (on peut aussi appliquer le principe de Lefschetz).
\end{remarque}
\section{Une propri\'et\'e de point fixe pour les fonctions rationnelles.}\label{propriete point fixe}
	Le but de cette section est de montrer la propri\'et\'e suivante.
	Rappelons que pour une fonction rationnelle $R \in \C_p(z)$ on note $R_*$ l'action sur $\H_p$ induite par $R$.
\begin{propriete de point fixe}
	Soit $R \in \C_p(z)$ une fonction rationnelle et soit $\hX \subset \H_p$ un ensemble connexe contenant au moins deux points, tel que $R_*(\hX) \subset \hX$.
	Alors, soit $\hX$ contient un point fixe {\rm rationnel} de $R_*$~; soit il existe un point fixe attractif $z_0 \in \pp$ de $R$, tel que $\hX$ contient une demi-g\'eod\'esique issue de $z_0$.
\end{propriete de point fixe}
	La d\'emonstration de cette propri\'et\'e et dans les Sections~\ref{lemmes preliminaires} et~\ref{sec propriete}.
	La Section~\ref{sec hyp} contient quelque rappels sur l'espace hyperbolique.

\

	On utilisera la Propri\'et\'e de Point Fixe de la fa\c{c}on suivante, voir Section~\ref{sec composantes periodiques}.
	Consid\'erons une fonction rationnelle $R \in \C_p(z)$ et un espace analytique $X \subset \pp$ tel que $R(X) \subset X$.
	Alors l'enveloppe convexe de $X$
$$
	\hX = \{ \cs \in \H_p \mid \cs \prec X \}
$$
(voir Section~\ref{enveloppes convexes}) est connexe et il contient au moins deux points (Lemme~\ref{detalles}).
	De plus la propri\'et\'e $R(X) \subset X$ implique $R_*(\hX) \subset \hX$ (Lemme~\ref{espaces extension}) et alors on peut appliquer la Propri\'et\'e de Point Fixe \`a $\hX$. 
\subsection{Quelque rappels sur l'espace hyperbolique.}\label{sec hyp}
	Dans cette section on consid\`ere de rappels sur l'espace hyperbolique $\H_p$, voir Sections~3 et~4 de~\cite{hyp} pour les r\'ef\'erences.

	L'espace hyperbolique $\H_p$ est munit d'un distance, qu'on note $d$, pour laquelle il est s\'eparable et complet.
	De plus $(\H_p, d)$ est un {\it arbre r\'eel}~: pour tous $\cs$, $\cs' \in \H_p$ distincts il existe un et un seul arc topologique ayant $\cs$ et $\cs'$ comme \'extremites et de plus cet arc topologique est isom\'etrique \`a un intervalle de $\R$.
	
	La distance entre deux points non singuliers $\cs$, $\cs' \in \H_p$ distincts est d\'efinit comme suit.
	Soit $B_\infty$ (resp. $B_\infty'$) la boule associ\'ee \`a $\cs$ (resp. $\cs'$) contenant $\infty$ et posons $D = \pp - B_\infty$ (resp. $D' = \pp - B_\infty$).
	Alors
\begin{equation}\label{formule distance}
	d(\cs, \cs') = \log_p \frac{\diam(D \cup D')^2}{\diam(D) \cdot \diam(D')},
\end{equation}
o\`u les diam\`etres sont par rapport \`a la distance sur $\C_p$ induite par la norme $| \cdot |$.
\subsubsection{Segments g\'eod\'esiques.}
\begin{definition}
	On dit que un point $\cs \in \H_p$ est {\bf entre} deux points distincts $\cs_0$ et $\cs_1$, si $\cs_0$ et $\cs_1$ appartienent \`a des composantes connexes distincts de $\H_p - \{ \cs \}$.

	On note $(\cs_0, \cs_1)$ l'ensemble de tous les points entre $\cs_0$ et $\cs_1$.
	De plus on pose $[\cs_0, \cs_1) = (\cs_1, \cs_0] = (\cs_0, \cs_1) \cup \{ \cs_0 \}$ et $[\cs_0, \cs_1] = [\cs_0, \cs_1) \cup \{ \cs_1 \}$.
\end{definition}
Pour $\cs_0$, $\cs_1 \in \H_p$ distincts l'ensemble $(\cs_0, \cs_1) \subset \H_p^\R$ est isom\'etrique \`a un intervalle de $\R$, on l'appelle {\it segment g\'eod\'esique}~; voir~\cite{hyp} Section~3.

	Si $\cs \in \H_p$ est un point singulier, alors l'ensemble $\H_p - \{ \cs \}$ est connexe et donc $\cs$ ne peut pas \^etre entre deux points.

	Si $\cs \in \H_p^\R$ est un point non singulier, alors les composantes connexes de $\H_p - \{ \cs \}$ sont les enveloppes convexes des boules associ\'ees \`a $\cs$ (voir Section~\ref{enveloppes convexes}).
	On a donc la partition
$$
	\H_p - \{ \cs \} = \sqcup_{\cs} \hB_\cp~;
$$
compare avec ($\ref{partition point}$) (rappelons que $\hB_\cp = \{ \cs \in \H_p \mid \cs \prec B_\cp \}$).

	Notons que si $\hX \subset \H_p$ est connexe, alors pour tous $\cs$, $\cs' \in \hX$ on a $(\cs, \cs') \subset \hX$.
	Par cons\'equent tout partie connexe de $\H_p$ contenant au moins deux points contient un point non singulier.

	Le lemme suivant c'est le Lemme~2.11 de~\cite{hyp}.
\begin{lemme}\label{direction Y}
	Soit $\cs \in \H_p$ et $\cs_0$, $\cs_1$ appartenant \`a la m\^eme composante connexe de $\H_p - \{ \cs \}$.
	Alors il existe un point $\ov\cs$ qui appartient \`a la m\^eme composante connexe, tel qu'on ait
$$
	[\cs_0, \cs) \cap [\cs_1, \cs) = [ \ov\cs, \cs).
$$
\end{lemme}
\subsubsection{Action des fonctions rationnelles sur $\H_p$.}
	Fixons une fonction rationnelle $R \in \C_p(z)$.
	Alors l'action $R_*$ sur $\H_p$ induite par $R$ satisfait
$$
	d(R_*(\cs), R_*(\cs')) \le \deg(R) \cdot d(\cs, \cs'), 
		\mbox{ pour } \cs, \cs' \in \H_p,
$$
voir Corollaire~4.7 de~\cite{hyp}.
	En particulier $R_*$ est continue.

	La proposition suivante suit de la Proposition~4.6 de~\cite{hyp}, voir aussi Lemme~\ref{direction Y}.
\begin{proposition}\label{boutdegre}
	Fixons un point $\cs \in \H_p$.
	\begin{enumerate}
		\item
			Pour tout $\cs_0 \in \H_P$ diff\'erent de $\cs$ il existe $\ov\cs \in (\cs, \cs_0)$ tel que $R_*$ soit injective sur $(\cs, \ov\cs)$ et tel qu'on ait $R_*((\cs, \ov\cs)) = (R_*(\cs), R_*(\ov\cs))$.
		\item
			Consid\'erons les param\'etrages
$$
	[\cs, \ov\cs] = \{ \cs(t) \}_{0 \ge t \ge r_0} \mbox{ et }
	[R_*(\cs), R_*(\ov\cs)] = \{ \ts(t) \}_{0 \ge t \ge \widetilde{r}_0}
$$
tel qu'on ait $d(\cs, \cs(t)) = t$ et $d(R_*(\cs), \ts(t)) = t$.
			Alors il existe un eniter $d \ge 1$ tel que pour $t > 0$ petit on ait $R_*(\cs(t)) = \ts(d \cdot t)$.
	\end{enumerate}
\end{proposition}
\subsection{Quelque lemmes pr\'eliminaires.}\label{lemmes preliminaires}
	La preuve du lemme suivant d\'epend de~\cite{hyp}.
\begin{lemme}\label{fixe => fixe rationnel}
	Soit $\hX \subset \H_p$ un ensemble connexe contenant au moins deux points et tel que $R_*(\hX) \subset \hX$.
	Si $R_*$ a un point fixe dans $\hX$, alors elle a aussi un point fixe {\rm rationnel} dans $\hX$.
\end{lemme}
\preuve
	Soit $\cs \in \hX$ un point fixe de $R_*$.
	On suppose que $\cs$ n'est pas rationnel.
	Comme $\hX$ contient au moins deux points il existe un point $\cs_0 \in \hX$ diff\'erent de $\cs$~; on a $(\cs, \cs_0) \subset \hX$, car $\hX$ est connexe.
	On montrera que tout point dans $(\cs, \cs_0) \subset \hX$, proche de $\cs$, est fix\'e par $R_*$.
	Comme les points rationnels sont denses sur chaque segment g\'eod\'esique (Lemme~3.3 de~\cite{hyp}) on obtient l'assertion du lemme.

	Si $\cs = \{ \cp \}$ est singulier alors $R_*(\cp) = \cp$ et par le Lemme~5.3 de~\cite{hyp} on a $\deg_R(\cp) = 1$.
	Alors par la partie 2 du Lemme~5.8 de~\cite{hyp} il existe un point $\ov\cs \in \H_p - \{ \cs \}$ tel que tout point de $(\cs, \ov\cs)$ soit fix\'e par $R_*$.
	Par le Lemme~\ref{direction Y} on peut supposer $\ov\cs \in (\cs_0, \cs)$.

	Si $\cs = \{ \cp, \cp' \}$ est un point irrationnel, alors $R_*(\cp) = \cp$, $R_*(\cp') = \cp'$ et $\deg_R(\cp) = \deg_R(\cp') = 1$ (Lemme~5.2 de~\cite{hyp}).
	On suppose $\cs_0 \prec B_\cp$.
	Alors l'assertion suit comme dans le cas pr\'ec\'edent.
$\hfill\square$
\begin{lemme}\label{po}
	Soient $\cs_0$, $\ts \in \H_p$ des points distincts et consid\'erons $\cs' \in (\cs_0, \ts)$.
	Alors pour tout point $\ts'$ satisfaisant $d(\ts, \ts') < d(\ts, \cs')$ on a $\cs' \in (\cs_0, \ts')$~; voir figure~\ref{figure1}.
\end{lemme}
\begin{figure}[htb]
\begin{center}
\psfig{file=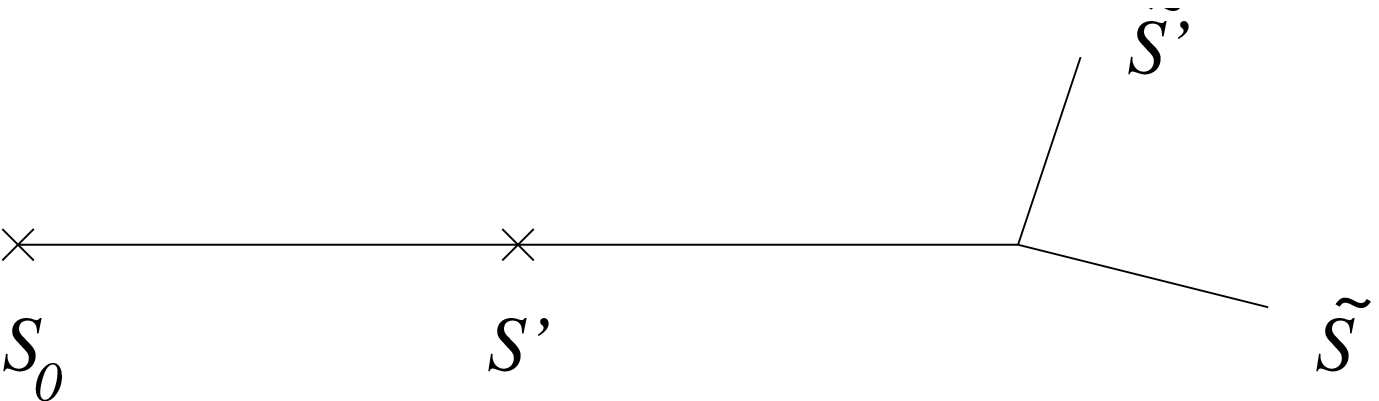, width = 3in}
\caption{}
\label{figure1}
\end{center}
\end{figure}
\preuve
	Il suffit de montrer que $\ts$ et $\ts'$ appartienent \`a la m\^eme composante connexe de $\H_p - \{ \cs \}$.
	Si ce ne pas le cas, alors $\cs' \in [\ts, \ts']$ et on aurait $d(\ts, \ts') \ge d(\ts, \cs')$.
$\hfill \square$
\begin{lemme}\label{demi-geodesique}
	Soit $\cs : [0, \infty) \longrightarrow \H_p$ une isom\'etrie sur son image, telle que $\cs(0) \in \H_p^\R$.
	Alors apr\`es changement de coordonn\'ee il existe $t_0 \in \R$ tel que pour $t \ge 0$, $\cs(t)$ est le point associ\'e \`a la boule $\{ |z| < p^{- t + t_0} \}$.
\end{lemme}
\preuve
	Notons que pour tout $t \ge 0$ le point $\cs(t)$ est non singulier.

	Apr\`es changement de coordonn\'ee on suppose que $\cs(0)$ est le point associ\'e \`a la boule $\{ |z| < p^{t_0} \}$ et que pour $t > 0$ le point $\cs(t)$ appartient \`a la composante connexe de $\H_p - \{ \cs(0) \}$ correspondante \`a la boule $\{ |z| < p^{t_0} \}$.
	C'est \`a dire $\cs(t) \prec \{ |z| < p^{t_0} \}$ pour $t > 0$.

	Pour $t \ge 0$ soit $B_\infty(t)$ la boule associ\'ee au point $\cs(t)$ qui contient $\infty$ et posons $D(t) = \pp - B_\infty(t)$.
	Comme $d(\cs(0), \cs(t)) = t$ on a $\diam(D(t)) = p^{-t + t_0}$~; cf. ($\ref{formule distance}$).
	En particulier $\diam(D(t))  \rightarrow 0$ lorsque $t \rightarrow \infty$.

	De plus pour $0 < t_0 < t_1$ le point $\cs(t_0)$ est entre $\cs(0)$ et $\cs(t_1)$, donc on a $D(t_1) \subset D(t_0)$.
	C'est \`a dire $\{ D(t) \}_{t \ge 0}$ est une suite decroissante de boules dont l'intersection $\cap_{t \ge 0} D(t)$ consiste d'un point (car $\C_p$ est complet).

	Apr\`es changement de coordonn\'ee fixant la boule $\{ |z| < p^{t_0} \}$, on suppose $\cap_{t \ge 0} D(t) = \{ 0 \}$.
	Alors $D(t) = \{ |z| \le p^{-t + t_0} \}$ et par cons\'equent $\cs(t)$ est le point associ\'e \`a la boule $\{ |z| < p^{-t + t_0} \}$.
$\hfill \square$
\subsection{Preuve de la Propri\'et\'e de Point Fixe.}\label{sec propriete}
	Fixons une fonction rationnelle $R \in \C_p(z)$ et soit $\hX \subset \H_p$ un ensemble connexe contenant au moins deux points et tel que $R_*(\hX) \subset \hX$.
	Supposons alors que $R_*$ n'a pas de point fixe rationnel dans $\hX$.
	Par le Lemme~\ref{fixe => fixe rationnel}, $R_*$ n'a pas de point fixe dans $\hX$.

\

$1.-$	On d\'efinira une application $\cs : [0, \infty) \longrightarrow \hX$ qui soit une isom\'etrie sur son image et telle pour chaque $t > 0$ on ait
$$
	\cs(t) \in (\cs(0), R_*(\cs(t))).
$$
	On choisit $\cs(0) \in \hX$ quelconque.
	Dans le 1.1 (resp. 1.2) on montrera que si l'application $\cs$ est d\'efinit pour $0 \le t \le t_0$ (resp. $0 \le t < t_0$), o\`u $t_0 \ge 0$, alors on peut \'etendre $\cs$ \`a $[0, t_1]$ (resp. $[0, t_0]$) o\`u $t_1 > t_0$.
	Ceci implique qu'on peut d\'efinir $\cs$ sur $[0, \infty)$.

\

$1.1.-$	Consid\'erons $t_0 \ge 0$ et supposons que l'application $\cs$ soit d\'ej\`a d\'efinit pour $0 \le t \le t_0$.
	On \'etendra $\cs$ \`a un intervalle $[0, t_1]$, o\`u $t_1 > t_0$.

	Comme $\cs(t_0) \in \hX$ on a $R_*(\cs(t_0)) \neq \cs(t_0)$ par hypoth\`ese.
	Soit $0 < \delta < d(R_*(\cs(t_0)), \cs(t_0))$ assez petit tel que pour tout point $\cs \in \H_p$ \`a distance au plus $\delta$ de $\cs(t_0)$ on ait
$$
	d(R_*(\cs(t_0)), R_*(\cs)) < d(R_*(\cs(t_0)), \cs(t_0)) - \delta.
$$
	Pour $0 \le t \le t_0 + \delta$ on d\'efinit $\cs(t)$ comme le seul point dans $(\cs(0), R_*(\cs(t_0)))$ \`a distance $t$ de $\cs(0)$.
	Alors l'application $\cs$ d\'efinit sur $[0, t_0 + \delta]$, est une isom\'etrie sur son image.

	Par d\'efinition de $\delta$, pour $t_0 < t \le t_0 + \delta$ on a
$$
	d(R_*(\cs(t_0)), R_*(\cs(t))) < d(R_*(\cs(t_0)), \cs(t_0)) - \delta
		\le d(R_*(\cs(t_0)), \cs(t)),
$$
or le Lemme~\ref{po} implique $\cs(t) \in (\cs(0), R_*(\cs(t)))$.

\

$1.2.-$	Supposons que l'application $\cs$ soit d\'ej\`a d\'efinit sur l'intervalle $[0, t_0)$, o\`u $t_0 > 0$.
	On \'etendra l'application $\cs$ \`a $[0, t_0]$.

	Comme $\cs : [0, t_0) \longrightarrow \H_p$ est une isom\'etrie sur son image et $\H_p$ est complet, il existe un point $\cs(t_0)$ limite de $\cs(t)$ lorsque $t \rightarrow t_0$.
	Comme il exsite un seul arc topologique dans $\H_p$ ayant $\cs(0)$ et $\cs(t_0)$ comme extr\'emit\'es (car $\H_p$ est un arbre r\'eel) l'application $\cs$ est un param\'etrage du segment $[\cs(0), \cs(t_0))$.

\

$1.2.1.-$	Supposons par contradiction $R_*(\cs(t_0)) = \cs(t_0)$.
	Par hypoth\`ese $R_*$ n'a pas de points fixes dans $\hX$, donc $\cs(t_0) \not\in \hX$.
	Soit $\hB$ la composante connexe de $\H_p - \{ \cs(t_0) \}$ contenant $\hX$.
	On a $\cs(0), R_*(\cs(0)) \in \hX \subset \hB$.
	Donc par la Proposition~\ref{boutdegre} et le Lemme~\ref{direction Y} il existe un entier $d \ge 1$ tel que pour $t < t_0$ proche de $t_0$ on ait
$$
	R_*(\cs(t)) = \cs(t_0 - d(t - t_0)).
$$
	Comme pour $t < t_0$ on a $t_0 - d(t - t_0) \le t$, ceci contredit la propri\'et\'e $\cs(t) \in (\cs(0), R_*(\cs(t)))$.

\

$1.2.2.-$	On a donc $R_*(\cs(t_0)) \neq \cs(t_0)$.
	On montrera d'abord que $\cs(t_0) \in (\cs(0), R_*(\cs(t_0)))$ et apr\`es on conclura $\cs(t_0) \in \hX$.

	Comme $R_*$ est continue on a
$$
	d(\cs(t), R_*(\cs(t))) \rightarrow d(\cs(t_0), R_*(\cs(t_0))) > 0
		\mbox{ lorsque } t \rightarrow t_0.
$$
	Donc pour $t < t_0$ proche de $t_0$ on a $d(R_*(\cs(t)), R_*(\cs(t_0))) < d(R_*(\cs(t)), \cs(t))$.
	Comme $\cs(t) \in (\cs(0), R_*(\cs(t)))$ le Lemme~\ref{po} implique dans ce cas $\cs(t) \in (\cs(0), R_*(\cs(t_0)))$.
	Alors
$$
	d(\cs(t_0), [\cs(0), R_*(\cs(t_0))]) \le d(\cs(t_0), \cs(t)) = t_0 - t.
$$
	Comme cette derni\`ere in\'egalit\'e est valable pour tout $t < t_0$ proche de $t_0$, on conclut $\cs(t_0) \in [\cs(0), R_*(\cs(t_0))]$.
	Comme $\cs(t_0) \neq \cs(0)$ et $\cs(t_0) \neq R_*(\cs(t_0))$, on a $\cs(t_0) \in (\cs(0), R_*(\cs(t_0)))$.

	Soit $t < t_0$ assez proche de $t_0$ tel que $d(R_*(\cs(t_0)), R_*(\cs(t))) < d(R_*(\cs(t_0)), \cs(t_0))$.
	Alors le Lemme~\ref{po} implique $\cs(t_0) \in (\cs(0), R_*(\cs(t)))$.
	Comme $\cs(0), \cs(t) \in \hX$ on a $R_*(\cs(t)) \in \hX$ et
$$
	\cs(t_0) \in (\cs(0), R_*(\cs(t))) \subset \hX.
$$

\

$2.-$	Comme $\hX \subset \H_p$ contient des points non singuliers (car $\hX$ est connexe et contient au moins deux points) on peut supposer $\cs(0)$ non singulier.
	Alors par le Lemme~\ref{demi-geodesique} apr\`es changement de coordonn\'ee le point $\cs(t)$ est le point associ\'e \`a la boule $\{ |z| < p^{- t + t_0} \}$.
	On montrera que $0$ est un point fixe attractif de $R$, ce qui termine la preuve.

\

$2.1.-$	Supposons par contradiction $R(0) \neq 0$.
	Apr\`es changement de coordonn\'ee on suppose $R(0) \in \C_p$.
	Notons que pour $t \gg 0$ l'ensemble $R(\{ |z| < p^{-t + t_0} \})$ est une boule, soit $r(t)$ son diam\`etre.
	Alors  $R(\{ |z| < p^{-t + t_0} \}) = \{ |z - R(0)| < r(t) \}$ est une boule associ\'ee au point $R_*(\cs(t))$.

	Si $t > 0$ est assez grand tel que $r(t) < |R(0)|$ alors le point $R_*(\cs(t))$ (qui est associ\'e \`a la boule $\{ |z - R(0)| < s(r) \}$) s\'epare la boule $\{ |z| > p^{-t + t_0} \} \cup \{ \infty \}$.
	Comme $\cs(0)$ aussi s\'epare la boule $\{ |z| > p^{-t + t_0} \} \cup \{ \infty \}$, les points $R_*(\cs(t))$ et $\cs(0)$ appartienent \`a la m\^eme composante connexe de $\H_p - \{ \cs(t) \}$.
	Mais ceci contredit la propri\'et\'e $\cs(t) \in (\cs(0), R_*(\cs(t)))$.

\

$2.2.-$	On conclut que $0 \in \C_p$ est un point fixe de $R$~; il reste \`a montrer qu'il est attractif.
	Posons $\lambda = R'(0)$.
	Si $\lambda = 0$ il n'y a rien \`a montrer, donc on suppose $\lambda \neq 0$.

	Alors pour $r > 0$ petit on a $R(\{ |z| < r \}) = R(\{ |z| < |\lambda| r \})$ et par cons\'equent pour $t \ge 0$ grand on a $R_*(\cs(t)) = \cs(t + \log_p|\lambda|)$.
	Donc la propri\'et\'e $\cs(t) \in (\cs(0), R_*(\cs(t)))$ implique $\log_p|\lambda| < 0$ et donc $|\lambda| < 1$.
	C'est \`a dire, $0$ est un point fixe attractif de $R$.
\section{Composantes p\'eriodiques de l'ensemble de Fatou.}\label{sec composantes periodiques}
	Fixons une fonction rationnelle $R \in \C_p(z)$ de degr\'e au moins deux.
	Dans cette section on deduit le Th\'eor\`eme~A~:
\begin{thmA}
	Soit $C_0 \subset \pp$ une composante p\'eriodique de $F(R)$.
	Alors il y a deux cas.
	\begin{enumerate}
		\item
			$C_0$ est un bassin d'attraction imm\'ediat.
		\item
			$C_0$ est une composante analytique du domaine de quasi-p\'eriodicit\'e.
	\end{enumerate}
\end{thmA}
	La preuve de ce th\'eor\`eme est dans la Section~\ref{preuve thmA}.

	Soit $\ca(R) \subset \pp$ l'union des bassins d'attraction imm\'eidat de $R$ et posons
$$
	\ca'(R) = \cup_{n \ge 0} R^{-n}(\ca(R)) \mbox{ et }
	\ce'(R) = \cup_{n \ge 0} R^{-n}(\ce(R)).
$$
	Notons que l'ensemble $\ca'(R)$ est l'union des bassins d'attraction de $R$.
	Il est facile de voir que les ensembles $\ca'(R)$ et $\ce'(R)$ sont disjoints.
	
	Notons que {\it le Th\'eor\`eme~A dit que $\ca(R) \sqcup \ce(R)$ (resp. $\ca'(R) \sqcup \ce'(R)$) est \'egal \`a l'union des composantes p\'eriodiques (resp. pr\'ep\'eriodiques) de l'ensemble de Fatou}.
	Le corollaire suivant est alors imm\'ediat.
\begin{corollaire}
	Pour une fonction rationnelle $R \in \C_p(z)$ les propri\'et\'es suivantes sont \'equivalentes.
	\begin{enumerate}
		\item
			$F(R) = \ca'(R) \sqcup \ce'(R)$.
		\item
			Il n'y a pas des composantes errantes de l'ensemble de Fatou.
	\end{enumerate}
\end{corollaire}
	On a conjectur\'e dans~\cite{these} qu'on a $F(R) = \ca'(R) \sqcup \ce'(R)$ pour toute fonction rationnelle.
	Par le corollaire pr\'ec\'edent ceci a lieu si et seulement si il n'y a pas des composantes errantes de l'ensemble de Fatou~; on peut comparer \`a~\cite{Be non rec}.

	Dans le cas complexe on sait, d'apr\`es un th\'eor\`eme de D.~Sullivan, que toute composante connexe de l'ensemble de Fatou est pr\'ep\'eriodique, voir~\cite{Su}.

	On a aussi le corollaire suivant du Th\'eor\`eme~A.
\begin{corollaire}
	Les composantes analytiques de $F(R) - (\ca(R) \sqcup \ce(R))$ co\"{\i}ndicent avec les composantes errantes de $F(R)$.
\end{corollaire}
\preuve
	Soit $C \subset \pp$ une composante errante de l'ensemble de Fatou.
	Comme $\ca'(R) \sqcup \ce'(R)$ c'est l'union des composantes pr\'ep\'eriodiques de l'ensemble de Fatou, on a $C \subset F(R) - \ca'(R) \sqcup \ce'(R)$.
	Par cons\'equent il existe une composante analytique $C'$ de cet ensemble qui contient $C$.
	On a
$$
	\ca'(C) \sqcup \ce'(R) \subset \pp - \cup_{n \ge 0} R^n(C').
$$
	Comme toute fonction rationnelle a un point fixe non r\'epulsif (voir~\cite{Be hyp}) on conclut que l'ouvert $\ca'(R) \sqcup \ce'(R)$ n'est pas vide et par cons\'equent il est de cardinalit\'e infinie.
	
	Par la d\'efinition de composante on a $C' \subset C$ et par cons\'equent $C' = C$.
$\hfill \square$
\subsection{Preuve du Th\'eor\`eme~A}\label{preuve thmA}
	Notons que par les Propositions~\ref{maximalite0} et~\ref{maximalite1} les bassins d'attraction imm\'ediat et les composantes analytiques du domaine de quasi-p\'eriodicit\'e sont des composantes de l'ensemble de Fatou.
	
	Alors le Th\'eor\`eme~A est une cons\'equence imm\'ediate de la proposition suivante (rappelons que l'ensemble de Fatou est l'union disjointe de ces composantes).
\begin{propA}
	Soit $X \subset \pp$ un espace anlytique tel que $R(X) \subset X$.
	Alors $X$ intersecte un bassin d'attraction imm\'ediat ou le domaine de quasi-p\'eriodicit\'e.
\end{propA}
	La preuve de cette proposition est plus bas, d'abord on consid\`ere quelque lemmes.
\begin{lemme}\label{disque}
	Soit $B = \{ |z| < r \}$ et soit $R \in \C_p(z)$ une fonction rationnelle telle que $R(B) \subset B$.
	\begin{enumerate}
		\item
			Si $R(B) = B$ et $R : B \longrightarrow R(B)$ est de degr\'e $1$, alors $B \subset {\cal E}(R)$ et pour tout $r_0 \in (0, r)$ et $r' > 0$ il existe $n \ge 1$ tel que $|R^n - id| \le r'$ sur $\{ |z| \le r_0 \}$.
		\item
			Si $R(B) \neq B$ ou $R : B \longrightarrow R(B)$ est de degr\'e plus grand que $1$, alors $B$ contient un point fixe attractif $z_0$ de $R$ et tout point dans $B$ converge vers $z_0$ par it\'eration positive.
	\end{enumerate}
\end{lemme}
{\par\noindent {\it Esquisse de preuve} ; voir~\cite{these}. \
	La partie $2$ est une cons\'equence simple du Lemme de Schwarz.
	Donc on suppose $R(B) = B$ et $R : B \longrightarrow B$ de degr\'e $1$.
	Notons que la deuxi\`eme assertion de la partie $1$ implique que tout point dans $B$ est r\'ecurrent par $R$ et donc $B \subset {\cal E}(R)$.
	
	Fixons $r_0 \in (|R(0)|, r)$ et consid\'erons la norme uniforme sur $B_0 = \{ |z| \le r_0 \}$, qu'on note $\| \cdot \|$.
	Comme $R : B \longrightarrow B$ est de degr\'e $1$ on a $|R'| = 1$ sur $B$ et il existe $\lambda \in \C_p$ de norme $1$, tel que $|R' - \lambda| < 1$ sur $B$.
	Quitte \`a changer $R$ par un it\'er\'e, on suppose $|R' - 1| < 1$ sur $B$.
	Alors il existe $\gamma \in (0, 1)$ tel que $\| R - id \| \le \gamma r$.
	
	Pour $m \ge 1$ on pose $\varepsilon_m = R^{p^m} - id$.
	Alors pour $k \ge 1$ on a,
$$
	R^{kp^m}(z) = z + \e_m(R^{p^m}) + \cdots + \e_m(R^{(k - 1)p^m})
		= k \e_m(z) + r_m(z),
$$
o\`u $\| r_m \| \le r^{-1} \| \e_m \|^2 \le \gamma \| \e_m \|$.
	Par cons\'equent on a $\| \e_{m + 1} \| = \frac{1}{p} \| \e_m \|$ pour $m$ assez grand.
	Donc $\| \e_m \| \longrightarrow 0$ lorsque $m \rightarrow \infty$.
$\hfill \square$
\begin{lemme}\label{fixe fatou}
	Soit $X \subset \pp$ un espace analytique.
	Supposons que $\cs \in \H_p^\Q$ est un point rationnel fix\'e par l'action de $R$ et tel que $\cs \prec X$.
	Alors $X$ intersecte un bassin d'attraction imm\'ediat ou le domaine de quasi-p\'eriodicit\'e.
\end{lemme}
\preuve
	Soit $\ct \subset \cs$ l'ensemble fini tel que pour tout $\cp \in \cs - \ct$ on ait $B_\cp \subset X$ et $R(B_\cp) = B_{R_*(\cp)}$ (voir Proposition~\ref{espace separation} et la partie 3 de la Proposition~\ref{point rationnel}).
	Comme toute fonction rationnelle \`a coefficients dans $\f$ \`a une infinit\'e de poins p\'eriodiques (cf. Lemme~10.1 de~\cite{hyp}) il existe un bout $\cp \in \cs$ p\'eriodique par $R_*$ et tel que $R_*^n(\cp) \not \in \ct$, pour $n \ge 0$.
	Alors il existe $k \ge 1$ tel que $R^k(B_\cp) = B_\cp$ or le lemme pr\'ec\'edent implique que $B_\cp \subset X$ intersecte un bassin d'attraction imm\'ediat ou le domaine de quasi-p\'eriodicit\'e.
$\hfill \square$

\

{\par\noindent{\it Preuve de la Proposition~A.}}
	Soit $\hX = \{ \cs \in \H_p^\R \mid \cs \prec X \}$ l'enveloppe convexe de $X$~; voir Sections~\ref{enveloppes convexes} et~\ref{propriete point fixe}.
	Alors $\hX$ est un ensemble connexe contenant au moins deux points (Lemme~\ref{detalles}).

	Comme $R(X) \subset X$, on a $R_*(\hX) \subset \hX$ (Lemme~\ref{espaces extension}).
	Alors, par la Propri\'et\'e de Point Fixe il y a deux cas.
	\begin{enumerate}
		\item
			Il existe un point fixe rationnel $\cs \in \H_p^\Q$ de $R_*$ tel que $\cs \prec X$.
			Alors la proposition suit du Lemme~\ref{fixe fatou}.
		\item
			Il existe un point fixe attractif $z_0$ de $R$, tel que $\hX$ contient une demi-g\'eod\'esique issue de $z_0$.
			Alors $z_0$ appartient \`a la fermeture topologique de $X$ (Lemme~\ref{detalles}) et par cons\'equent $X$ intersecte le bassin d'attraction imm\'ediat de $z_0$.
	\end{enumerate}
$\hfill \square$
\section{Appendice 1.}
	L'objectif de cet appendice est de montrer le lemme suivante, qu'on utilise dans la Section~\ref{sec repulsif}.
\begin{lemme}\label{couronne1}
	Consid\'erons $r_0, r_0' > 0$ et posons $B_0 = \{ |z| < r_0 \}$ et $B_0' = \{ |z| < r_0' \}$.
	Supposons que l'image par l'action de $R$ du bout de $B_0$ est le bout de $B_0'$.
	Consid\'erons $r \in (0, r_0)$ et soit $C = \{ r < |z| < r_0 \}$.
	Alors on a exactement une des propri\'et\'es suivantes.
	\begin{enumerate}
		\item
			$B_0' \subset R(C)$.
		\item
			L'ensemble $C' = R(C) \cap B_0'$ est une couronne dont l'un des ces bouts est le bout de $B_0'$.

			Si de plus le degr\'e de $R$ au bout de $B_0$ est au moins deux, alors il existe $m > 0$ ind\'ependent de $r$, tel que on ait
$$
	\mod(C') \ge \mod(C) + \min \{ m, \mod (C) \}.
$$
	\end{enumerate}
\end{lemme}
	Le preuve de ce lemme d\'epend du lemme suivant.
\begin{lemme}\label{couronne0}
	Soient $0 < r_1 < r_0$ et supposons que la fonction rationnelle $R$ n'a pas ni z\'eros ni p\^oles dans $\{ r_1 < |z| < r_0 \}$.
	Alors $|R(z)| = a|z|^d$ sur $\{ r_1 < |z| < r_1 \}$, o\`u $a > 0$ et $d \in \Z$.
	Si l'on a $d \ge 1$, alors $R(\{ r_1 < |z| < r_0 \}) = \{ ar_1^d < |z| < a r_0^d \}$.
\end{lemme}
\preuve
	Par hypoth\`ese $R$ est de la forme $\lambda P/Q$, o\`u $\lambda \in \C_p - \{ 0 \}$,
$$
	P(z) = (z - a_1) \cdot \ldots \cdot (z - a_k) \cdot
		(1 - z b_1) \cdot \ldots \cdot (1 - zb_m),
$$
$$
	Q(z) = (z - a_1') \cdot \ldots \cdot (z - a_{k'}') \cdot
		(1 - z b_1') \cdot \ldots \cdot (1 - zb_{m'}'),
$$
avec $|a_i|$, $|a_i'| < r_0$ et $|b_j|$, $|b_j'| < r_1^{-1}$.
	Alors $r_0 < |z| < r_1$ implique $|R(z)| = |\lambda| \cdot |z|^{k - k'}$.
	
	Supposons $d = k - k' \ge 1$.
	Apr\`es changement de coordonn\'ee au d\'epart on peut supposer $|\lambda| = |w|$.
	Alros tout $r_0 < |z| < r_1$ satisfaisant $R(z) = w$ est de norme \'egale \`a $1$.
	De plus $P(z) - w Q(z) = \lambda S(z)$, o\`u $S(z) = z^k - (w/\lambda) z^{k'} \mod \, {\frak m}_p$ (rappelons que ${\frak m}_p$ est l'ideal $\{ |z| < 1 \}$).
	Par cons\'equent il y a $d = k - k' \ge 1$ pr\'eimages de $w$ dans $r_0 < |z| < r_1$, compt\'ees avec multiplicit\'e.
$\hfill \square$

\

{\par\noindent{\it Preuve du Lemme~\ref{couronne1}.}}
	Pour $r > 0$ soit $\cs_r$ le point associ\'e \`a $\{ |z| < r \}$.
	Supposons $B_0' \not\subset R(C)$~; apr\`es changement de coordonn\'ee on se ram\`ene au cas $0 \not \in R(C)$.
	Soient $r < r_k < \ldots < r_1 < r_0$ les normes des p\^oles de $R$ dans $C$.
	On pose $r = r_{k + 1}$ et pour $0 \le i \le k$ on pose $C_i = \{ r_{i + 1} < |z| < r_i \}$.

\

$1.-$	Comme $R$ n'a pas ni z\'eros ni p\^oles dans $C_i$, on a $|R(z)| = a_i |z|^{d_i}$ sur $C_i$, o\`u $a_i > 0$ et $d_i \in \Z$ (Lemme~\ref{couronne0}).
	Comme l'image du bout de $B_0$ est le bout de $B_0'$, on a $d_0 \ge 1$.
	Donc $R(C_0) = \{ r_1' < |z| < r_0' \}$, o\`u $r_1' = a_0r_1^{d_0}$ (Lemme~\ref{couronne0}).

\

$2.-$	On montre par induction $R(C_i) = \{ r_{i + 1}' < |z| < r_i' \}$ pour $0 \le i \le k$, o\`u $0 < r_{k + 1}' < \ldots < r_1' < r_0'$.
	Par le Lemme~\ref{couronne0} il suffit de montrer $d_i \ge 1$ pour $0 \le i \le k$.
	
	Par $1$, ceci est vraie pour $i = 0$~; supposons par induction que cette propri\'et\'e est vraie pour $i - 1$.
	Alors l'image du bout de $\{ |z| > r_i \} \cup \{ \infty \}$ est le bout de $\{ |z| > r_i' \} \cup \{ \infty \}$ et par cons\'equent $R_*(\cs_{r_i}) = \cs_{r_i'}$.
	Soit $B$ une boule associ\'ee \`a $\cs_{r_i}$ telle que l'image de son bout soit le bout de $\{ |z| < r_i' \}$.
	On a $0 \in \{ |z| < r_i' \} \subset R(B)$ (Lemme~\ref{residu}) et comme $0 \not \in R(\{ r < |z| < r_0 \})$ on conclut $B = \{ |z| < r_i \}$.
	C'est \`a dire, l'image du bout de $\{ |z| < r_i \}$ est le bout de $\{ |z| < r_i' \}$ et par cons\'equent $d_i \ge 1$.

\

$3.-$	Soit $1 \le i \le k$ et soit $B \subset \{ |z| = r_i \}$ une boule associ\'ee \`a $\cs_{r_i}$.
	Comme $0 \not \in R(C)$, $R(B)$ est une boule (Lemme~\ref{residu}) qui est associ\'ee \`a $\cs_{r_i'}$ et qui est diff\'erente de $\{ |z| < r_i' \}$.
	Comme $R_* : \cs_{r_i} \longrightarrow \cs_{r_i'}$ est surjective (Proposition~\ref{point rationnel}), on conclut que $R(\{ |z| = r_i \})$ est soit \'egal \`a $\{ |z| = r_i' \}$, soit \'egal \`a $\{ |z| \ge r_i' \} \cup \{ \infty \}$.
	Dans les deux cas on a
$$
	C' = R(C) \cap B_0' = R(\{ r_{k + 1} < |z| < r_0 \}) \cap \{ |z| < r_0' \} =
		\{ r_{k + 1}' < |z| < r_0' \}.
$$
	
\

$4.-$	Il reste \`a montrer la derni\`ere assertion de la partie 2.
	Comme l'image du bout de $B_0$ est le bout de $B_0'$, il existe $s \in (0, r_0)$ ind\'ependent de $r$ tel que on ait $|R(z)| = a_0 |z|^{d_0}$ sur $\{ s < |z| < r_0 \}$.
	Soit $m = \mod( \{ s < |z| < r_0 \})$, on a $\mod(C_0) \ge \min\{ m, \mod(C) \}$.
	Notons que $d_0$ est le degr\'e de $R$ au bout de $B_0$.
	Donc par $2$, si $d_0 > 1$ on a
\begin{eqnarray*}
	 \mod(C) &  = & \sum d_i \cdot \mod(C_i) \ge \sum \mod(C_i) + (d_0 - 1)\cdot \mod(C_0). \\
		& \ge & \mod(C) + \min \{ m, \mod(C) \}.
\end{eqnarray*}
$\hfill \square$

\bibliographystyle{plain}

\end{document}